\PassOptionsToPackage{hypertexnames=false}{hyperref}
\documentclass[]{Liebert_Author_revised}

\newtheorem{hypothesis}[theorem]{\bf{Hypothesis}}
\newtheorem{definition}{Definition}
\usepackage{capt-of}
\usepackage{indentfirst}
\usepackage{xcolor}

\title{From the SIR/SEIR/SVEIR Model to a Chemotactic System: Analysis of Spatiotemporal Mechanisms of Infection Spread within the Einstein Paradigm} 

\author{
Yizhou~Wang,$^{1\ast}$ Evgenia~Echkina,$^{1}$\\
{$^{1}$Department of Automation for Scientific Research,}\\
{Faculty of Computational Mathematics and Cybernetics,}\\
{Lomonosov Moscow State University,}\\
{Leninskie Gory 1-52, Moscow 119991, Russia}\\
{$^\ast$To whom correspondence should be addressed;}\\
{E-mail: saipeizyu@gmail.com.}\\
{Co-author e-mail: ejanester@gmail.com.}
}

\date{}
\begin{document}
\singlespacing

\maketitle 


\begin{abstract}
This work investigates the class of SIR/SEIR models and their extension to an ODE-based SVEIR system, whose parameters are calibrated and validated using real influenza outbreak data.

The theoretical part establishes existence and uniqueness of solutions for the SVEIR model via the Picard--Lindel\"of theorem and analyzes invariance of the admissible region.

To describe spatiotemporal infection spread in an urban setting, a diffusion--chemotaxis extension of the SIR framework is formulated as a PDE system; an additional term is introduced to represent fear/avoidance effects in the system dynamics.

A numerical workflow is proposed, including grid discretization, Neumann boundary conditions enforced by mirror extension, a semi-implicit IMEX time-stepping scheme, and sparse linear solves; computational correctness is monitored, in particular, by checking conservation of the total population.
\end{abstract}

\section{Introduction and Research Background}

Classical compartmental epidemic models continue to provide a fundamental mathematical framework for studying the transmission of infectious diseases. The historical development of this framework is commonly traced to Hamer's mass-action idea in 1906 and Ross's malaria model in 1911, before the work of Kermack and McKendrick placed the susceptible--infected--recovered structure on a systematic mathematical footing \citep{Brauer_2017,Kermack_McKendrick_1927}.

Subsequent developments incorporated additional biological mechanisms. In particular, the introduction of an exposed class led to the SEIR framework; the earliest explicit SEIR formulation is commonly traced to Cooke's 1965/1967 work \citep{Burke_2024}. Later, vaccine-related epidemic modeling motivated the SVEIR class, including SARS-oriented formulations that explicitly incorporate vaccination effects \citep{Gumel_McCluskey_Watmough_MBE_2006,Gao_Huang_Kang_Zhang_BVP_2018}.

In this paper, we revisit the SIR, SEIR, and SVEIR frameworks and derive an ODE-based SVEIR system under a short-term assumption in which vaccination is not treated as a dominant driving factor. We then consider epidemic propagation in an urban environment and formulate a constant-diffusion SIR system of PDEs together with a fear-augmented SIR PDE model. These models are used to investigate spatiotemporal transmission dynamics under population mobility.

The spatial part of the study is also related to diffusion-based descriptions of movement and to more recent work connecting Brownian-motion ideas with chemotactic and nonlinear PDE systems. Einstein's 1905 Brownian-motion paper provided a differential-equation-based description of microscopic random motion \citep{Einstein_1905}. More recently, Ibragimov and collaborators have used Einstein's method to motivate chemotactic systems, while related work has examined propagation-speed issues in degenerate Einstein-type diffusion models \citep{Islam_Ibragimov_CMFD_2024,GarliHevage_Ibraguimov_Sobol_arXiv_2023}. Building on this line of ideas, the present paper further develops an urban epidemic diffusion framework for the susceptible, infected, and recovered populations.

\section{Epidemiological Models with ODEs}
\subsection{SIR Model}

The \textbf{SIR} model is a fundamental conceptual framework in epidemiological modeling, which make our model similar to "prey-predator" paradigm. Its core idea is to partition the population into three mutually exclusive groups according to an individual's infection status: \textbf{susceptible} ($S$), \textbf{infectious} ($I$), and \textbf{recovered} ($R$). The model is built on the assumption that the total population size remains constant throughout the epidemic, thereby neglecting demographic processes such as births, deaths, and migration. The key mechanism driving the model dynamics is the transmission of the pathogen from infectious individuals to susceptible individuals.

The dynamics of the SIR model are described by a system of ordinary differential equations. Let the total population size be $N$, and let the numbers of individuals in each compartment at time $t$ be denoted by $S(t)$, $I(t)$, and $R(t)$, respectively. 

The transitions between these compartments are governed by
\begin{equation}
\tag{1.1}\label{eq:1.1}
\begin{aligned}
\frac{dS}{dt} &= -\beta\,\frac{S I}{N},\\
\frac{dI}{dt} &= \beta\,\frac{S I}{N}-\gamma I,\\
\frac{dR}{dt} &= \gamma I,\\
N &= S+I+R.
\end{aligned}
\end{equation}

The parameter $\beta$ denotes the contact intensity, which determines the probability of infection transmission per effective contact. The parameter $\gamma$, in turn, represents the mean recovery rate; its reciprocal, $1/\gamma$, provides an estimate of the average duration over which an infected individual remains capable of transmitting the disease.

We will perform numerical simulations using computer code. To enhance the practical relevance of the model demonstration, we will estimate (fit) the model parameters based on a real dataset from an influenza outbreak in a British boarding school, as reported in the \emph{British Medical Journal} in 1978.

The daily case data used for fitting are
\texttt{[1, 3, 6, 25, 73, 222, 294, 258, 237, 191, 125, 69, 27, 11, 4]}.
The dataset was taken from \citep{Musacchio_SIRModel_Blog_2020}.

\subsection{Solution Method for the SIR Model (RK4)}\label{sec:sir-rk4}

When solving the SIR system, instead of using standard built-in numerical routines in Python---such as the library function \texttt{odeint} adopted in \cite{Musacchio_SIRModel_Blog_2020}---we implemented the classical \textbf{fourth-order Runge--Kutta (RK4)} method to numerically integrate this system of ordinary differential equations. The purpose of this choice is to explicitly present the step-by-step computational procedure and to prepare a transparent implementation framework for subsequently solving our proposed multivariate models.

Here the state vector is defined as $u=[S,I,R]^\top$, and the right-hand side of the system is given by
$f(t,u)=\left[-\frac{\beta SI}{N},\ \frac{\beta SI}{N}-\gamma I,\ \gamma I\right]^\top$.

\textbf{Step 1: Computation of $k_{1}$:}\quad
\begin{equation}
\begin{aligned}
k_{1} =h\,f(t_{n},u_{n}) 
= h\left[
-\frac{\beta S_{n}I_{n}}{N},\ 
\frac{\beta S_{n}I_{n}}{N}-\gamma I_{n},\ 
\gamma I_{n}
\right]^\top.
\end{aligned}
\tag{2.1}
\end{equation}

\textbf{Step 2: Computation of $k_{2}$ and $k_{3}$:}

Intermediate state:
\begin{equation}
\begin{split}
u_{\mathrm{mid}1} =u_{n}+\frac{k_{1}}{2}
=\left[
S_{n}+\frac{k_{1,S}}{2},\ 
I_{n}+\frac{k_{1,I}}{2},\ 
R_{n}+\frac{k_{1,R}}{2}
\right]^\top.
\end{split}
\tag{2.2}
\end{equation}
\begin{equation}
\begin{split}
u_{\mathrm{mid}2} =u_{n}+\frac{k_{2}}{2}
=\left[
S_{n}+\frac{k_{2,S}}{2},\ 
I_{n}+\frac{k_{2,I}}{2},\ 
R_{n}+\frac{k_{2,R}}{2}
\right]^\top.
\end{split}
\tag{2.3}
\end{equation}

\textbf{Then:}\quad
\begin{equation}
\begin{aligned}
&k_{2} = h\,f\!\left(t_{n}+\frac{h}{2},\,u_{\mathrm{mid}1}\right) 
      = h\left[
-\frac{\beta\left(S_{n}+\frac{(k_{1})_{S}}{2}\right)\left(I_{n}+\frac{(k_{1})_{I}}{2}\right)}{N},\right. \\
&\left.
\frac{\beta\left(S_{n}+\frac{(k_{1})_{S}}{2}\right)\left(I_{n}+\frac{(k_{1})_{I}}{2}\right)}{N} \right. 
\left.
-\gamma\left(I_{n}+\frac{(k_{1})_{I}}{2}\right),
\gamma\left(I_{n}+\frac{(k_{1})_{I}}{2}\right)
\right]^{\top}.
\end{aligned}
\tag{2.4}
\end{equation}

\begin{equation}
\begin{aligned}
k_{3} &= h\,f\!\left(t_{n}+\frac{h}{2},\,u_{\mathrm{mid}2}\right) 
     = h\left[
-\frac{\beta\left(S_{n}+\frac{(k_{2})_{S}}{2}\right)\left(I_{n}+\frac{(k_{2})_{I}}{2}\right)}{N}, \right. \\
&\left.
\frac{\beta\left(S_{n}+\frac{(k_{2})_{S}}{2}\right)\left(I_{n}+\frac{(k_{2})_{I}}{2}\right)}{N} \right.
\left.
-\gamma\left(I_{n}+\frac{(k_{2})_{I}}{2}\right), 
\gamma\left(I_{n}+\frac{(k_{2})_{I}}{2}\right)
\right]^{\top}.
\end{aligned}
\tag{2.5}
\end{equation}

\textbf{Step 3: Computation of $k_{4}$:}\quad

Intermediate state:
\begin{equation}
\begin{aligned}
u_{\mathrm{end}}=u_{n}+k_{3}
=\left[S_{n}+(k_{3})_{S},\ I_{n}+(k_{3})_{I},\ R_{n}+(k_{3})_{R}\right]^{\top}.
\end{aligned}
\tag{2.6}
\end{equation}

\textbf{Then:}\quad
\begin{equation}
\begin{aligned}
& k_{4} =h\,f(t_{n}+h,\,u_{\mathrm{end}}) 
 =h\left[
-\frac{\beta\left(S_{n}+(k_{3})_{S}\right)\left(I_{n}+(k_{3})_{I}\right)}{N},\right. \\
&\left.
\frac{\beta\left(S_{n}+(k_{3})_{S}\right)\left(I_{n}+(k_{3})_{I}\right)}{N} \right. 
\left.
-\gamma\left(I_{n}+(k_{3})_{I}\right),
\gamma\left(I_{n}+(k_{3})_{I}\right)
\right]^{\top}.
\end{aligned}
\tag{2.7}
\end{equation}

\textbf{Step 4: State update:}\quad
\begin{equation}
\tag{2.8}
u_{n+1}=u_{n}+\frac{1}{6}\left(k_{1}+k_{2}+k_{3}+k_{4}\right).
\end{equation}

\textbf{Componentwise form.}
\begin{equation}
\tag{2.9}
\begin{aligned}
S_{n+1} &= S_{n}+\frac{1}{6}\left(k_{1,S}+k_{2,S}+k_{3,S}+k_{4,S}\right),\\
I_{n+1} &= I_{n}+\frac{1}{6}\left(k_{1,I}+k_{2,I}+k_{3,I}+k_{4,I}\right),\\
R_{n+1} &= R_{n}+\frac{1}{6}\left(k_{1,R}+k_{2,R}+k_{3,R}+k_{4,R}\right).
\end{aligned}
\end{equation}

Figure~\ref{fig:sir_ode} illustrates the numerical trajectories of the SIR ordinary differential equation model obtained using the above Runge--Kutta scheme.
The analysis shows that, for a fixed infection transmission coefficient $\beta = 0.4$, the model predictions exhibit a significant deviation from the empirical data with respect to the peak incidence, which reaches 294 cases on day 36.
Reducing the transmission coefficient to $\beta = 0.3$ results in a noticeably improved agreement between the model output and the observed data.
However, such a parameter adjustment contradicts the epidemiological context of the corresponding period, since the limited development of medical infrastructure would reasonably imply a higher effective transmission intensity.

\noindent\begin{minipage}{\linewidth}
    \centering
    \includegraphics[width=\linewidth]{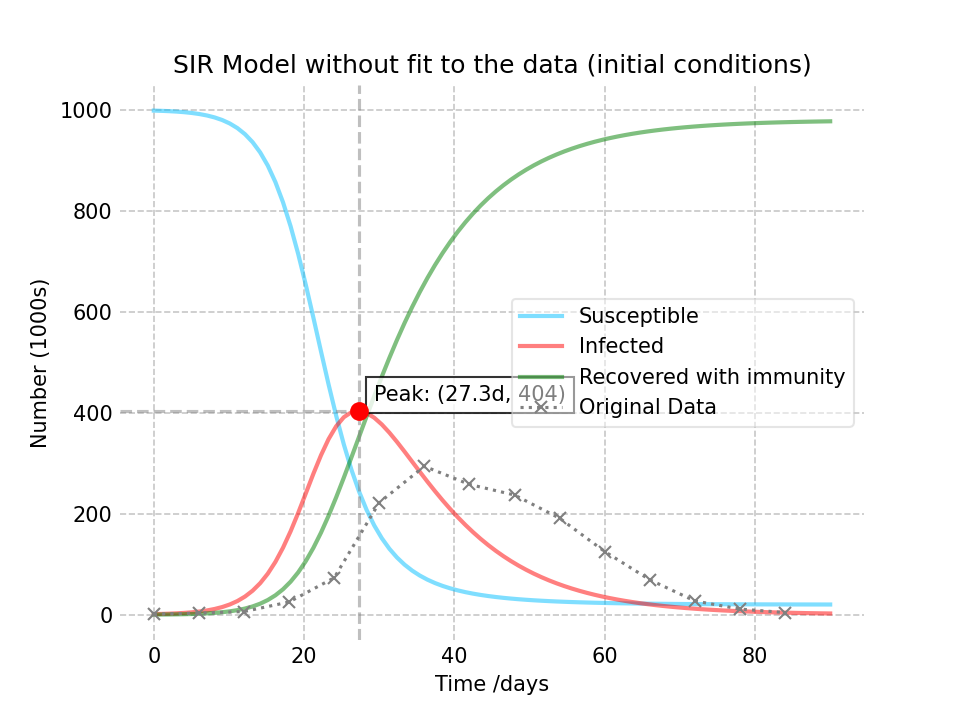}
    \captionof{figure}{SIR trajectories of an ordinary differential equation (ODE) model under the initial (non-fitted) parameter setting.}
    \label{fig:sir_ode}
\end{minipage}

\subsection{SEIR Model}

As is well known, the classical SIR model is well established and has long been employed in scientific practice.
In recent years, an extended SEIR model has been widely adopted in the epidemiological literature.
The key modification consists in introducing an additional population compartment, denoted by $E$ (exposed individuals), into the classical SIR framework.

This compartment represents infected individuals who are in a latent stage of the disease and have not yet developed clinical symptoms.
The inclusion of the exposed class allows the model to capture the delay between infection and infectiousness and, without altering the basic parameter set, to achieve a more accurate representation of the observed epidemic dynamics.
The structure of the extended SEIR model is given by
\begin{equation}
\tag{1.2}\label{eq:1.2}
\begin{aligned}
\frac{dS}{dt} &= -\beta \frac{SI}{N},\\
\frac{dE}{dt} &= \beta \frac{SI}{N} - \sigma E,\\
\frac{dI}{dt} &= \sigma E - \gamma I,\\
\frac{dR}{dt} &= \gamma I,\\
N &= S + E + I + R .
\end{aligned}
\end{equation}

For the numerical simulation of the SEIR system, the same initial dataset and baseline parameter values as those used for the SIR model were retained in order to ensure a consistent comparison between the two formulations. This makes it possible to isolate the effect of introducing the exposed compartment without introducing additional variability through changes in the empirical input or in the principal epidemiological parameters. The average duration of the latent period was assumed to be two days, corresponding to the choice $\sigma^{-1}=2$.

The resulting system of ordinary differential equations was solved numerically by means of the classical fourth-order Runge--Kutta method, following the same computational procedure as described in Subsection~\ref{sec:sir-rk4}. This provides a transparent numerical framework for tracking the temporal evolution of all four compartments and enables a direct comparison with the corresponding SIR dynamics.

The numerical results show that the SEIR model exhibits a closer agreement with the observed outbreak data than the basic SIR model. From an epidemiological perspective, this improvement is explained by the presence of the exposed class, which accounts for the latent stage of infection and thereby allows the model to represent the delay between infection and infectiousness more realistically. As a consequence, the overall temporal profile of the epidemic is reproduced more accurately.

Figure~\ref{fig:seir_ode} displays the numerical trajectories of the SEIR ordinary differential equation model under the initial (non-fitted) parameter setting. The figure shows the temporal evolution of the susceptible, exposed, infected, and recovered populations and illustrates the improved qualitative correspondence of the SEIR framework with the observed epidemic dynamics.

\noindent\begin{minipage}{\linewidth}
    \centering
    \includegraphics[width=\linewidth]{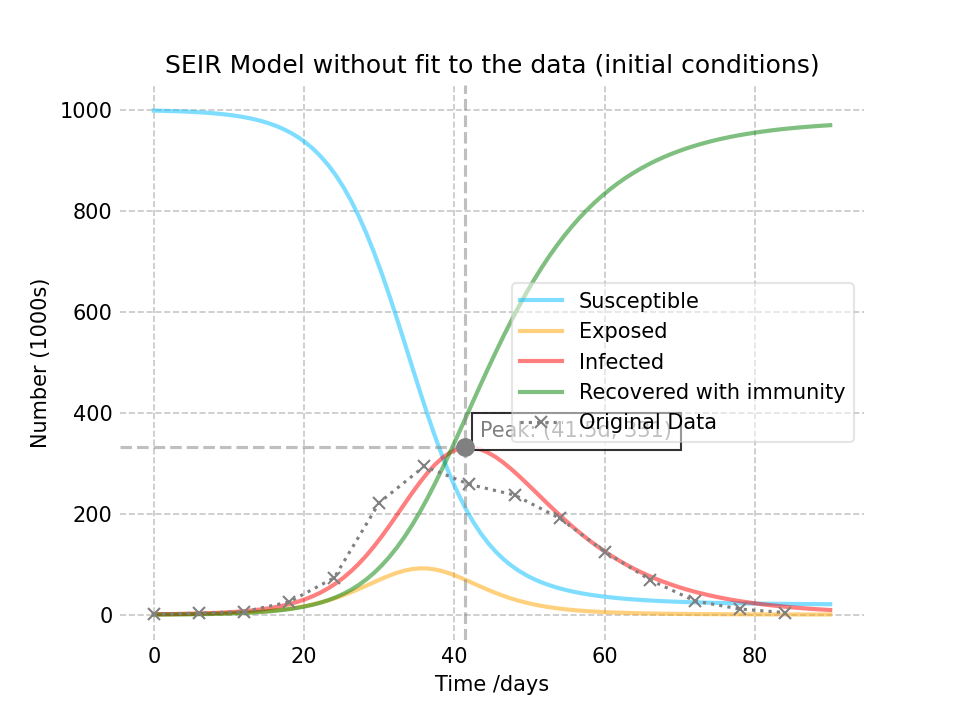}
    \captionof{figure}{SEIR trajectories of an ordinary differential equation (ODE) model under the initial (non-fitted) parameter setting.}
    \label{fig:seir_ode}
\end{minipage}

\subsection{The SVEIR model and its construction}

\subsubsection{The SVEIR model}

With the advancement of medical science, an important question naturally arises: which additional factors may influence the evolution of the number of infected individuals?
If one assumes that a fraction of the population has been vaccinated, this may lead not only to positive effects but also to complex and difficult-to-predict consequences.
To account for this factor, we extend Model~(2) by introducing a new population compartment $V$ (vaccinated individuals) and further develop the model accordingly.

In the scientific literature, a wide variety of models incorporating a vaccinated compartment ($V$) have been proposed, in which diffusion-based approaches are commonly employed. In contrast, we aim to modify the baseline SEIR model by directly introducing additional differential equations to construct an SVEIR model. This approach provides a systematic foundation for the subsequent development of an epidemiological model with extended structure and greater flexibility. After a systematic analysis, the following system of equations is obtained.

\begin{equation}
\tag{1.3}\label{eq:1.3}
\begin{aligned}
\frac{dS}{dt} &= -\beta \frac{SI}{N} - \nu S,\\
\frac{dV}{dt} &= \nu S - (1-\varepsilon)\beta \frac{VI}{N},\\
\frac{dE}{dt} &= \beta \frac{SI}{N} + (1-\varepsilon)\beta \frac{VI}{N} - \sigma E,\\
\frac{dI}{dt} &= \sigma E - \gamma I,\\
\frac{dR}{dt} &= \gamma I,\\
N &= S + V + E + I + R .
\end{aligned}
\end{equation}

Here, the parameter $\nu$ denotes the vaccination rate, while $\varepsilon$ represents the vaccine efficacy.
However, since the dataset employed in this study originates from the 1970s, a period during which medical infrastructure and vaccination coverage were relatively limited, both parameters are assumed to take comparatively low values.

The model is solved in a manner analogous to the two previous models.
After a minor adjustment of the numerical code, we observe that the model outcomes differ depending on the size of the vaccinated population. In the simulations presented below, a very low vaccination rate $\nu = 0.004$ is selected, and the vaccine efficacy is assumed to be $\varepsilon = 0.6$ (60\%). So we get:

\noindent\begin{minipage}{\linewidth}
    \centering
    \includegraphics[width=\linewidth]{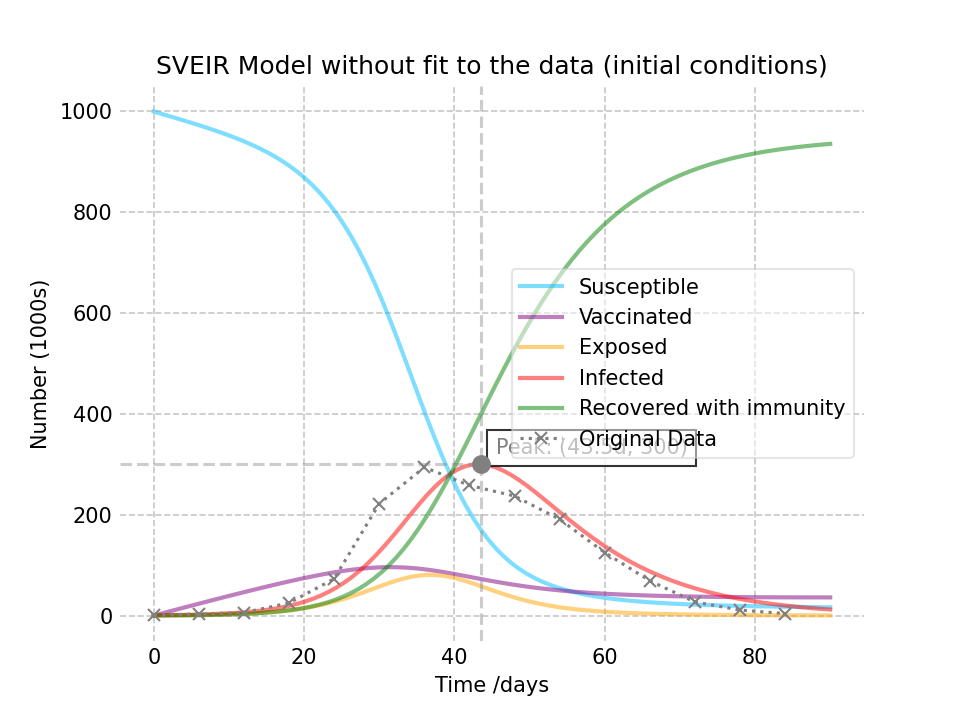}
    \captionof{figure}{SVEIR trajectories of an ordinary differential equation (ODE) model under the initial (non-fitted) parameter setting.}
    \label{fig:sveir_ode}
\end{minipage}

In Figure~\ref{fig:sveir_ode}, we saw that, according to the model defined by model~\eqref{eq:1.3}, constructed using the given dataset, the estimated peak incidence is reached approximately on day~43.5 with about 300 cases, which is consistent with the epidemiological situation of that period.
In contrast, for Model~(2), which does not incorporate vaccination, the peak incidence occurs on day~41.5 with 331 cases.
Based on a comparison of these two sets of results, it can be inferred that vaccination against the virus was most likely not implemented in that British hospital at the time.
With regard to contemporary medical practice, we argue that the SVEIR model is more relevant, as increasing emphasis is placed on preventive medicine.

\subsubsection{From the Baseline SVEIR ODE Formulation to a One-Dimensional Spatiotemporal Model via the Einstein Paradigm}

We first introduce the Einstein mathematical paradigm for $1$ dimension spacial variable, as follows. 
\begin{equation}
\begin{aligned}
\vec{u}(x,t+\tau)
= \vec{u}(x,t)
+ \tau\,\vec{F}\bigl(\vec{u}(x,t)\bigr)
+ \int_{\mathbb{R}^d}
\bigl[
\vec{u}(x+\zeta,t)-\vec{u}(x,t)
\bigr]\,
\varphi(\zeta,x,t)\,d\zeta .
\end{aligned}
\tag{2.10}
\end{equation}

In Model~(1.3), the vector $\vec{u}$ is a column vector composed of five components,
\begin{equation}
\tag{2.11}
\vec{u}=(S,V,E,I,R)^{\top}.
\end{equation}

By expanding the integrand, we obtain
\begin{equation}
\begin{aligned}
&\vec{u}(x,t+\tau) - \vec{u}(x,t)
\;+\; \int_{\mathbb{R}^d}
\vec{u}(x,t)\,\varphi(\zeta,x,t)\,d\zeta 
= \tau\,\vec{F}\bigl(\vec{u}(x,t)\bigr)
+ \int_{\mathbb{R}^d}
\vec{u}(x+\zeta,t)\,\varphi(\zeta,x,t)\,d\zeta .
\end{aligned}
\tag{2.12}
\end{equation}

That is,
\begin{equation}
\begin{aligned}
&\vec{u}(x,t+\tau) - \vec{u}(x,t)
\;+\; \vec{u}(x,t)
\int_{\mathbb{R}^d}
\varphi(\zeta,x,t)\,d\zeta
= \tau\,\vec{F}\bigl(\vec{u}(x,t)\bigr)
+ \int_{\mathbb{R}^d}
\vec{u}(x+\zeta,t)\,\varphi(\zeta,x,t)\,d\zeta .
\end{aligned}
\tag{2.13}
\end{equation}

Here, $\varphi$ represents a probability density (jump kernel) and satisfies the following conditions:
\begin{equation}
\tag{2.14}
\int_{\mathbb{R}^d} \varphi(\zeta)\,d\zeta = 1,
\qquad
0 \le \varphi(\zeta) \le 1 .
\end{equation}

At this stage, the left-hand side of the previous equation can be transformed as
\begin{equation}
\begin{aligned}
&\vec{u}(x,t+\tau) - \vec{u}(x,t) \;
+\; \vec{u}(x,t)\int_{\mathbb{R}^d}\varphi(\zeta,x,t)\,d\zeta
= \vec{u}(x,t+\tau) - \vec{u}(x,t) + \vec{u}(x,t)\cdot 1
= \vec{u}(x,t+\tau).
\end{aligned}
\tag{2.15}
\end{equation}

Consequently, the left-hand side becomes a partial derivative with respect to time $t$.
Observing that the operation of expanding the integral is redundant and does not contribute to computational tractability, we refrain from performing this expansion.
In accordance with Model~(1.3), which does not involve spatial variables and depends solely on time, the integral of the difference term on the right-hand side can be taken to be identically zero.

Dividing both sides of the equation by the parameter $\tau$, the right-hand side reduces to a single vector-valued function $\vec{F}(\vec{u})$, that is,
\begin{equation}
\begin{aligned}
\frac{\vec{u}(t+\tau)-\vec{u}(t)}{\tau}
\;\xrightarrow[\tau \to 0]{}\;
\frac{\partial \vec{u}}{\partial t}(t)
= \vec{F}\bigl(\vec{u}(t)\bigr).
\end{aligned}
\tag{2.16}
\end{equation}

We may now employ this expression to describe the first five equations of Model~(1.3).

\subsubsection{Existence and Uniqueness of Solutions for the SVEIR Model (Theorem 2.1)}
\textbf{Theorem 2.1.(Existence and uniqueness of solutions for the SVEIR model)}

Consider the SVEIR model written as the system of ordinary differential equations~(3).
Let the initial data be prescribed at time $t_0\ge 0$ by
\begin{equation}
\begin{aligned}
S(t_0)=S_0,\quad V(t_0)=V_0,\quad E(t_0)=E_0, \quad
I(t_0)=I_0,\quad R(t_0)=R_0.
\end{aligned}
\tag{3.1}
\end{equation}
Assume that
\begin{equation}
\begin{aligned}
S_0\ge 0,\quad V_0\ge 0,\quad E_0\ge 0,\quad
I_0\ge 0,\quad R_0\ge 0,
\end{aligned}
\tag{3.2}
\end{equation}
and
\begin{equation}
\begin{aligned}
S_0+V_0+E_0+I_0+R_0 = N>0.
\end{aligned}
\tag{3.3}
\end{equation}
Then the following statements hold:
\begin{enumerate}
\item \textbf{(Local existence and uniqueness).}
For the given initial configuration $(S_0,V_0,E_0,I_0,R_0)$, there exists a unique solution
\[
\bigl(S(t),V(t),E(t),I(t),R(t)\bigr)
\]
of system~(3), defined at least on  time interval $[t_0,t_0+\tau]$ for some $\tau>0$.

\item \textbf{(Global continuation and invariance of the feasible region).}
The above solution can be uniquely continued to the entire half-line $t\ge t_0$.
Moreover, for all $t\ge t_0$ the following relations hold:
\begin{equation}
\begin{aligned}
S(t)\ge 0,\quad V(t)\ge 0,\quad E(t)\ge 0,\quad
I(t)\ge 0,\quad R(t)\ge 0,
\end{aligned}
\tag{3.4.1}
\end{equation}
and
\begin{equation}
\begin{aligned}
S(t)+V(t)+E(t)+I(t)+R(t)=N>0.
\end{aligned}
\tag{3.5.1}
\end{equation}
\end{enumerate}

\begin{proof} 
\textbf{The Proof of Theorem 2.1.}

To prove this theorem, we invoke a classical result from the theory of ordinary differential equations, namely the Picard--Lindel\"of theorem on existence and uniqueness for the Cauchy (initial value) problem. \\
“\textbf{(Picard--Lindel\"of theorem on existence and uniqueness).}

Consider the Cauchy problem
\begin{equation}
\begin{aligned}
u'(t)=f\bigl(t,u(t)\bigr), \qquad u(t_0)=u_0,
\end{aligned}
\tag{3.6}
\end{equation}
where $u(t)\in\mathbb{R}^n$. Suppose there exists a rectangular domain
\begin{equation}
\begin{aligned}
Q=\{(t,u): |t-t_0|\le T,\ \|u-u_0\|\le A\},
\end{aligned}
\tag{3.7}
\end{equation}
in which the function $f(t,u)\in C(Q)$ is continuous in $(t,u)$ and satisfies the Lipschitz condition with respect to $u$, that is, there exists a constant $L>0$ such that for all $(t_1,u_1),(t_2,u_2)\in Q$,
\begin{equation}
\begin{aligned}
\|f(t,u_1)-f(t,u_2)\|\le L\,\|u_1-u_2\|.
\end{aligned}
\tag{3.8}
\end{equation}
Then there exists $h\in(0,T]$ such that on the interval $[t_0-h,\,t_0+h]$ the Cauchy problem admits a unique solution $u(t)$.
Moreover, the solution can be obtained as the limit of the Picard sequence:
\begin{equation}
\begin{aligned}
u_0(t)=u_0,
\end{aligned}
\tag{3.9}
\end{equation}
\begin{equation}
\begin{aligned}
u_{n+1}(t)=u_0+\int_{t_0}^{t} f\bigl(\tau,u_n(\tau)\bigr)\,d\tau,
\end{aligned}
\tag{3.10}
\end{equation}
and $u_n\to u$ uniformly on $[t_0-h,\,t_0+h]$.”

First, we introduce the vector
\begin{equation}
\begin{aligned}
u(t)=\bigl(S(t),V(t),E(t),I(t),R(t)\bigr)^{T}\in\mathbb{R}^{5}.
\end{aligned}
\tag{3.11}
\end{equation}
Then the system can be written compactly as the Cauchy problem
\begin{equation}
\begin{aligned}
\frac{du}{dt}=F(u),\qquad u(0)=u_0,
\end{aligned}
\tag{3.12}
\end{equation}
where $u_0=(S_0,V_0,E_0,I_0,R_0)^{T}$ is the initial data, and $F:\mathbb{R}^{5}\to\mathbb{R}^{5}$ is defined componentwise.

where
\begin{equation}
\begin{aligned}
F(u)=
\begin{pmatrix}
f_{S}(S,V,E,I,R)\\
f_{V}(S,V,E,I,R)\\
f_{E}(S,V,E,I,R)\\
f_{I}(S,V,E,I,R)\\
f_{R}(S,V,E,I,R)
\end{pmatrix}.
\end{aligned}
\tag{3.13.1}
\end{equation}
For instance,
\begin{align*}
f_{S}(S,V,E,I,R) &= -\beta\,\frac{SI}{N}-\nu S, \\
f_{I}(S,V,E,I,R) &= \sigma E-\gamma I.
\end{align*}
and similarly for $f_V$, $f_E$, and $f_R$.

Now consider any component, for instance $f_S$ for fixed $S$. It is a linear combination of products and sums: $SI$, $S$, and the constants $\beta$, $\nu$, and $N$. Such expressions are polynomials (or rational functions with the constant denominator $N>0$) and, therefore, are continuous with respect to all five variables.

The same is true for the remaining components $f_V$, $f_E$, $f_I$, and $f_R$. Each of these expressions is obtained from $S$, $V$, $E$, $I$, and $R$ by means of addition, subtraction, and multiplication by constants. Hence, each component $f_S$, $f_V$, $f_E$, $f_I$, and $f_R$ is continuous, and consequently the vector field
\begin{equation}
\begin{aligned}
F(u)=(f_S,f_V,f_E,f_I,f_R)
\end{aligned}
\tag{3.13.2}
\end{equation}
is continuous on $\mathbb{R}^5$.

The Picard--Lindel\"of theorem requires a local Lipschitz condition with respect to $u$. One convenient way to verify this is to note that each component $f_i$ has continuous partial derivatives with respect to $S$, $V$, $E$, $I$, and $R$.

For example, for the component $f_S$ the partial derivatives are
\[
\frac{\partial f_S}{\partial S}=-\beta\,\frac{I}{N}-\nu,\qquad
\frac{\partial f_S}{\partial I}=-\beta\,\frac{S}{N},
\]
while with respect to the remaining variables $(V,E,R)$ the derivatives are equal to zero. Similarly, one can write down the derivatives for $f_V$, $f_E$, $f_I$, and $f_R$. All of them are linear combinations of $S$, $V$, $E$, and $I$ with constants, for instance,
\[
\frac{\partial f_V}{\partial S}=\nu,\qquad
\frac{\partial f_V}{\partial V}=(\varepsilon-1)\beta\,\frac{I}{N},\qquad
\frac{\partial f_V}{\partial I}=(\varepsilon-1)\beta\,\frac{V}{N},
\]
and so on.

Now let us fix a bounded domain in which the solution will be considered. A natural choice for an epidemiological model is
\begin{equation}
\begin{aligned}
\mathcal{Q}=\left\{u=(S,V,E,I,R)\in\mathbb{R}^{5}:\ 0\le S,V,E,I,R\le N\right\}.
\end{aligned}
\tag{3.14}
\end{equation}

On this domain all variables are bounded. Hence, all partial derivatives
$\dfrac{\partial f_i}{\partial S}$, $\dfrac{\partial f_i}{\partial V}$, $\dfrac{\partial f_i}{\partial E}$, \dots
are bounded by some common constant $B>0$.
It is known from multivariable analysis that if a function $F:\mathbb{R}^{5}\to\mathbb{R}^{5}$ is continuously differentiable and its Jacobian is bounded on the domain $\mathcal{Q}$, then $F$ satisfies the Lipschitz condition on $\mathcal{Q}$; that is, there exists a constant $L>0$ such that
\begin{equation}
\begin{aligned}
\|F(u_1)-F(u_2)\|\le L\,\|u_1-u_2\|,\quad
\text{for all } u_1,u_2\in\mathcal{Q}.
\end{aligned}
\tag{3.15}
\end{equation}

Consequently, the function $F(u)$ defining the right-hand side of the SVEIR model system satisfies the local Lipschitz condition with respect to $u$.
By this theorem, there exists $h>0$ such that on the time interval $[0,h]$ there exists a unique solution $u(t)$ to our Cauchy problem. In other words, the trajectory
\begin{equation}
\begin{aligned}
t\mapsto \bigl(S(t),V(t),E(t),I(t),R(t)\bigr)
\end{aligned}
\tag{3.16}
\end{equation}
exists and is unique at least in some neighborhood of the initial time.

Now we establish two simple but important properties of the solutions.

\textbf{1. The sum of the components is conserved.}
Adding all equations, we obtain
\begin{equation}
\begin{aligned}
&\frac{d}{dt}\bigl(S+V+E+I+R\bigr)
=\left(-\beta\frac{SI}{N}-\nu S\right)
+\left(\nu S-(1-\varepsilon)\beta\frac{VI}{N}\right) \\
&\quad+\left(\beta\frac{SI}{N}+(1-\varepsilon)\beta\frac{VI}{N}-\sigma E\right) +(\sigma E-\gamma I) +(\gamma I).
\end{aligned}
\tag{3.17}
\end{equation}
All terms cancel pairwise, and we obtain
\begin{equation}
\begin{aligned}
\frac{d}{dt}\bigl(S+V+E+I+R\bigr)=0.
\end{aligned}
\tag{3.18}
\end{equation}
Therefore,
\begin{equation}
\begin{aligned}
S(t)+V(t)+E(t)+I(t)+R(t)\equiv N,
\end{aligned}
\tag{3.19}
\end{equation}
for all $t$ for which the solution exists.

\textbf{2. Nonnegativity of the components (invariance of the positive orthant).}
Intuitively, all flows between the compartments are described by terms of the form ``constant $\times$ current value''; when a variable reaches zero, the corresponding ``outflow'' vanishes.

For example, if at some time $t_1$ we have $S(t_1)=0$, then
\begin{equation}
\begin{aligned}
\frac{dS}{dt}(t_1)=-\beta\,\frac{0\cdot I(t_1)}{N}-\nu\cdot 0=0,
\end{aligned}
\tag{3.20}
\end{equation}
that is, starting from the state $S=0$ the quantity cannot ``drop'' into the negative region. Analogous arguments apply to $V$, $E$, $I$, and $R$ (taking into account the balance of incoming and outgoing flows). A more rigorous formulation is as follows: the set
\begin{equation}
\begin{aligned}
\mathcal{Q}_{+}
=\{u=(S,V,E,I,R)\in\mathbb{R}^{5}: S,V,E,I,R\ge 0, 
\left. \ \ S+V+E+I+R=N\right\}
\end{aligned}
\tag{3.21}
\end{equation}
is invariant for the system; that is, if $u(0)\in\mathcal{Q}_{+}$, then $u(t)\in\mathcal{Q}_{+}$ for all $t$ in the interval of existence of the solution.

From these two facts it follows that the solution trajectory does not leave the bounded region:
\begin{equation}
\begin{aligned}
0\le S(t),V(t),E(t),I(t),R(t)\le N,
\end{aligned}
\tag{3.4.2}
\end{equation}
\begin{equation}
\begin{aligned}
S(t)+V(t)+E(t)+I(t)+R(t)=N.
\end{aligned}
\tag{3.5.2}
\end{equation}

The Picard--Lindel\"of theorem yields a local solution on $[0,h]$. However, we are interested in the epidemic dynamics on the entire interval $t\ge 0$. A general principle of ODE theory states that if the right-hand side is Lipschitz on the region where the solution evolves and the solution remains bounded, then the solution can be continued arbitrarily far in time (that is, there is no ``blow-up'' in finite time).

In our case:
\begin{itemize}
\item we already know that the solution remains in $\mathcal{Q}_{+}$ for all times, and this is a bounded region (all coordinates lie between $0$ and $N$);
\item on the set $\mathcal{Q}_{+}$ the function $F(u)$ remains Lipschitz (the gradients are bounded).
\end{itemize}

Consequently, the solution does not ``leave'' the region on which the right-hand side is well defined, and we may apply the existence and uniqueness theorem step by step, extending the interval of existence. As a result, the local solution on $[0,h]$ can be continued to any finite interval $[0,T]$, and hence to the entire half-line $t\ge 0$. Therefore, for model~\eqref{eq:1.3} global existence and uniqueness hold: for every admissible initial state $u_0\in\mathcal{Q}_{+}$ there exists a unique trajectory $u(t)$ defined for all $t\ge 0$, which preserves the nonnegativity of the components and the total population size $N$.

\end{proof}

For model~\eqref{eq:1.3} (SVEIR), the right-hand side of the system is continuous and locally Lipschitz with respect to the state variables; therefore, by the Picard--Lindel\"of theorem, the initial value problem admits a unique local solution. Using the invariance of the nonnegative region and the conservation of the total population size $N$, we conclude that the solution remains bounded and can be continued on the entire interval $t\ge 0$. Thus, model~\eqref{eq:1.3} is well posed: for any biologically meaningful initial data, there exists a unique solution of the system that is global in time.

\subsection{Comparison of Results for Different Models}\label{subsec:comparison-models}

For convenience in comparing the SEIR model~\eqref{eq:1.2} and the SVEIR model~\eqref{eq:1.3}, we set the initial number of vaccinated individuals equal to zero, which reflects the absence of a developed vaccine against the new infectious disease at that time. We compare several groups corresponding to different vaccination levels: at the zero vaccination level we recover the standard SEIR model~\eqref{eq:1.2}, whereas all subsequent cases correspond to the SVEIR model~\eqref{eq:1.3}. It is evident that the higher the vaccination level, the flatter the incidence curve becomes. However, given the real conditions of that period, one cannot assert with certainty whether a corresponding vaccine had been developed; nevertheless, the construction of the SVEIR model~\eqref{eq:1.3} and the numerical experiments already provide a solid foundation for future studies of infectious diseases.

\noindent\begin{minipage}{\linewidth}
\centering
\includegraphics[width=\linewidth]{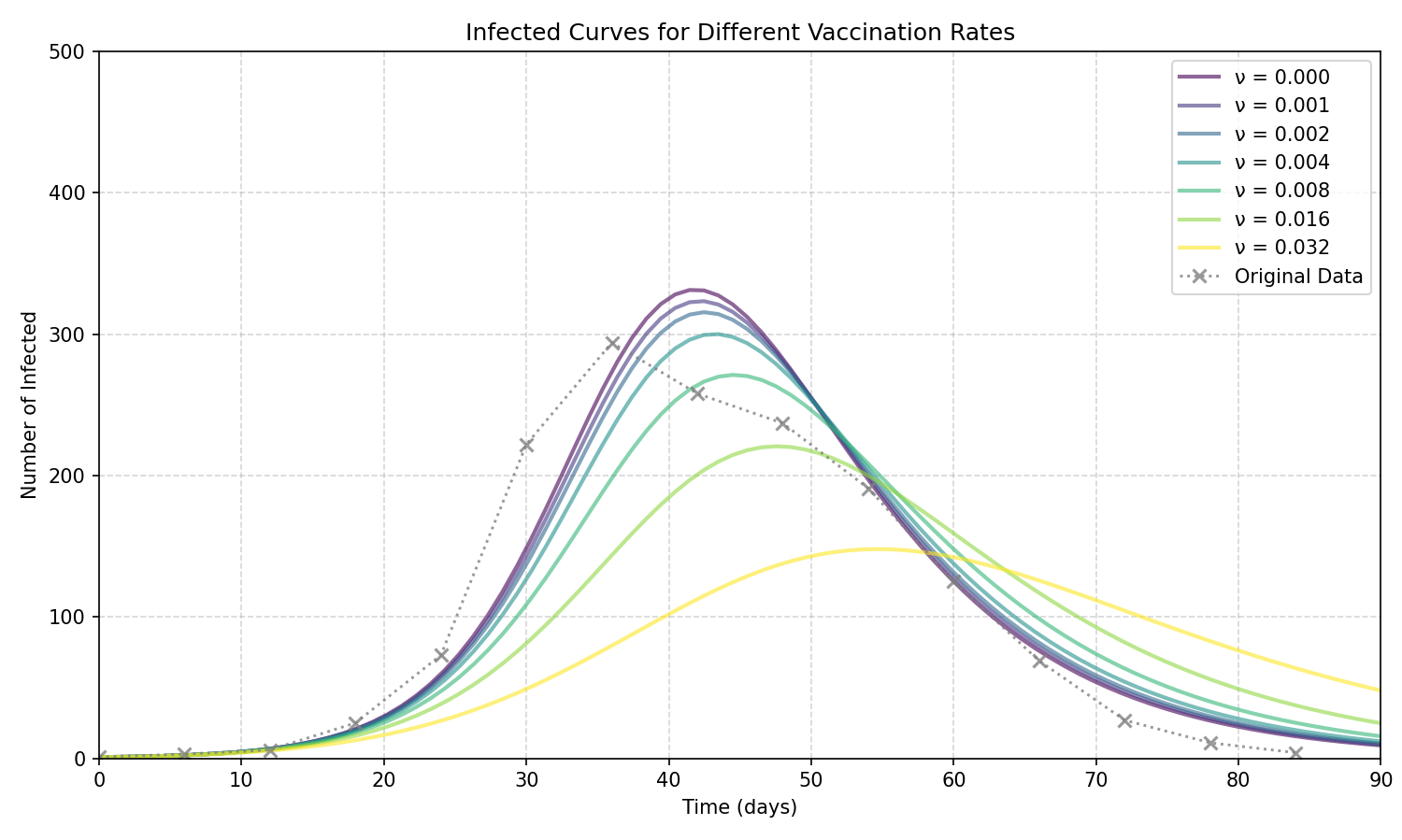}
\captionof{figure}{Infected curves for different vaccination rates (comparison of the SEIR model~\eqref{eq:1.2} and the SVEIR model~\eqref{eq:1.3}).}
\label{fig:compare-12-13}
\end{minipage}

\section{Spatial Dynamics in the SIR Model: Accounting for Diffusion}\label{sec:spatial-sir-diffusion}

\subsection{\texorpdfstring{Construction of a One-Dimensional SIR Model with a Chemotactic Term via the Einstein Paradigm (Model~\eqref{eq:1.6.1}--\eqref{eq:1.7.1})}{Construction of a One-Dimensional SIR Model with a Chemotactic Term via the Einstein Paradigm (Model 1.6.1--1.7.1)}}

Then we propose that 
\begin{equation}
\boxed{\begin{aligned}
&S(x,t+\tau) - S(x,t)\
\;
= \; S(x,t)\frac{\partial }{\partial x}\int_{\mathbb{R}^1}
\zeta\varphi_S(\zeta,x,t)\,d\zeta +\tau_S\,\vec{F_S}\bigl(S,I(x,t)\bigr)\\
&+ \int_{\mathbb{R}^1}
(S(x+\zeta,t)-S(x,t))\varphi_{S}(\zeta,x,t)\,d\zeta .
\end{aligned}
}
\tag{4.1}
\end{equation}

Let us hypothesize that my chemotactic force to be such that 
\begin{hypothesis}
\begin{equation}
\boxed{\frac{\partial }{\partial x}\int_{\mathbb{R}^1}
\zeta\varphi_S(\zeta,x,t)\,d\zeta=-\chi\frac{\frac{\partial I}{\partial x}}{I}}
\tag{4.2}
\end{equation}.
\end{hypothesis}

Let us introduce condition for all probability functions of $S, I, R$.
Namely we assume that all functions $\varphi$ are such that,
\begingroup
\renewcommand{\thedefinition}{2.1}
\begin{definition}\label{def:phi_condition}
\[
\ \varphi(\zeta,x,t) \ \text{ to be such that}
\ \varphi(\zeta,x,t)=0  \  \text{for} \ \zeta<0
\]
\end{definition}
\endgroup

Taking the second equation and definition \ref{def:phi_condition},  our SIR transport model takes the following form form
\begin{equation}  
\boxed{\begin{aligned}
&I(x,t+\tau) - I(x,t) 
= \tau_I\,\vec{F_I}\bigl(S,I(x,t)\bigr)\\
&+ \int_{\mathbb{R}^1}
(I(x+\zeta,t)-I(x,t))\varphi_{I}(\zeta,x,t)\,d\zeta ,
\end{aligned}
}
\tag{4.3}
\end{equation}
and finally equation for recovered species will take a form
\begin{equation}
\boxed{\begin{aligned}
&R(x,t+\tau) - R(x,t)
= \tau_R\,\vec{F_R}\bigl(R(x,t),I(x,t)\bigr)\\
&+ \int_{\mathbb{R}^1}
(R(x+\zeta,t)-R(x,t))\varphi_{R}(\zeta,x,t)\,d\zeta .
\end{aligned}
}
\tag{4.4}
\end{equation}

Not that recovered and infected species move only in diffusive fashion, namely we assume that functions $\varphi_I$, and $\varphi_R$ are independent on $I$.

As a consequence, we obtain a first--order SIR model with a chemotactic term:
\begin{equation}
\begin{alignedat}{1}
&\tau_S\,\frac{\partial S}{\partial t}
=-\beta\,\frac{SI}{N}
+D_S\,\frac{\partial^2 S}{\partial x^2}
+\chi\,\frac{\partial}{\partial x}\!\left(S\,\frac{1}{I}\,\frac{\partial I}{\partial x}\right),\\
&\tau_I\,\frac{\partial I}{\partial t}
=\beta\,\frac{SI}{N}
-\gamma I
+D_I\,\frac{\partial^2 I}{\partial x^2},\\
&\tau_R\,\frac{\partial R}{\partial t}
=\gamma I
+D_R\,\frac{\partial^2 R}{\partial x^2},\\
&S(x,t)+I(x,t)+R(x,t)=N(x,t)>0.
\end{alignedat}
\tag{1.6.1}\label{eq:1.6.1}
\end{equation}
Its initial and boundary conditions are specified as follows. Using homogeneous Neumann (zero--flux) conditions, we obtain:
\begin{equation}
\begin{alignedat}{1}
&x\in[0,L]:\quad S_x(0,t)=S_x(L,t)=0,\quad I_x(0,t)=I_x(L,t)=0,\quad R_x(0,t)=R_x(L,t)=0,\\
&S(x,0)=A e^{-B(x-x_0)^2},\qquad I(x,0)=E e^{-F(x-x_1)^2},\qquad R(x,0)=0.
\end{alignedat}
\tag{1.7.1}\label{eq:1.7.1}
\end{equation}

\begingroup
\newcounter{defsave}
\setcounter{defsave}{\value{definition}} 

\renewcommand{\thedefinition}{2.\arabic{definition}}
\setcounter{definition}{0} 


\setcounter{definition}{\value{defsave}} 
\endgroup

\subsection{SIR Reaction--Diffusion System without a Chemotactic Term}\label{subsec:sir-no-chemotaxis}
\subsubsection{\texorpdfstring{Model~\eqref{eq:1.4}--\eqref{eq:1.5}}{Model 1.4--1.5}}\label{subsubsec:sir-no-chemotaxis-model}

Now we shift the analysis to alternative scenarios. The previous data on the influenza outbreak in a British boarding school were based on the assumption of a closed population, which was justified under those conditions. However, when studying communities or cities, population mobility must be taken into account. How can the spread of infection be tracked in such settings? We may introduce a diffusion term to model human movement and incorporate two parameters characterizing the spatial boundaries. For simplicity, we begin this extension from the baseline model~\eqref{eq:1.1}, obtaining:

\begin{equation}
\begin{aligned}
\frac{\partial S}{\partial t} &= -\beta\,\frac{S\,I}{N} + D_S\,\Delta S,\\
\frac{\partial I}{\partial t} &= \beta\,\frac{S\,I}{N} - \gamma\,I + D_I\,\Delta I,\\
\frac{\partial R}{\partial t} &= \gamma\,I + D_R\,\Delta R,\\
S+I+R &= N>0,\\
\Delta(\cdot) &= \frac{\partial^2(\cdot)}{\partial x^2}+\frac{\partial^2(\cdot)}{\partial y^2}.
\end{aligned}
\tag{1.4}\label{eq:1.4}
\end{equation}

For the boundary conditions, we impose homogeneous Neumann conditions and specify the initial conditions:
\begin{equation}
\begin{alignedat}{1}
&\left.\frac{\partial S}{\partial n}\right|_{\Gamma}
=\left.\frac{\partial I}{\partial n}\right|_{\Gamma}
=\left.\frac{\partial R}{\partial n}\right|_{\Gamma}
=0,\\[2pt]
&S(x,y,0)=S_0(x,y) 
=A\exp\!\Bigl[-B\bigl((x-x_0)^2
+(y-y_0)^2\bigr)\Bigr],\\[2pt]
&I(x,y,0)=1,\\
&R(x,y,0)=0,\\[2pt]
&\frac{\partial}{\partial n}
=\cos(\vec n,x)\frac{\partial}{\partial x}
+\cos(\vec n,y)\frac{\partial}{\partial y}.
\end{alignedat}
\tag{1.5}\label{eq:1.5}
\end{equation}

And as a result, the total population size is conserved, which is confirmed by the numerical results:
\begin{equation}
\begin{aligned}
\mathrm{Total}=\iint S\,dx\,dy.
\end{aligned}
\tag{2.17}
\end{equation}
Since a Gaussian distribution was adopted for the initial conditions, the total population size satisfies
\begin{equation}
\begin{aligned}
\mathrm{Total}\approx \frac{A\pi}{B}.
\end{aligned}
\tag{2.18}
\end{equation}

Under the proposed model and the adopted dataset/parameter setting, the total population remains (approximately) conserved over time. The corresponding evolution of the total population is shown in:

\noindent\begin{minipage}{\linewidth}
\centering
\includegraphics[width=\columnwidth]{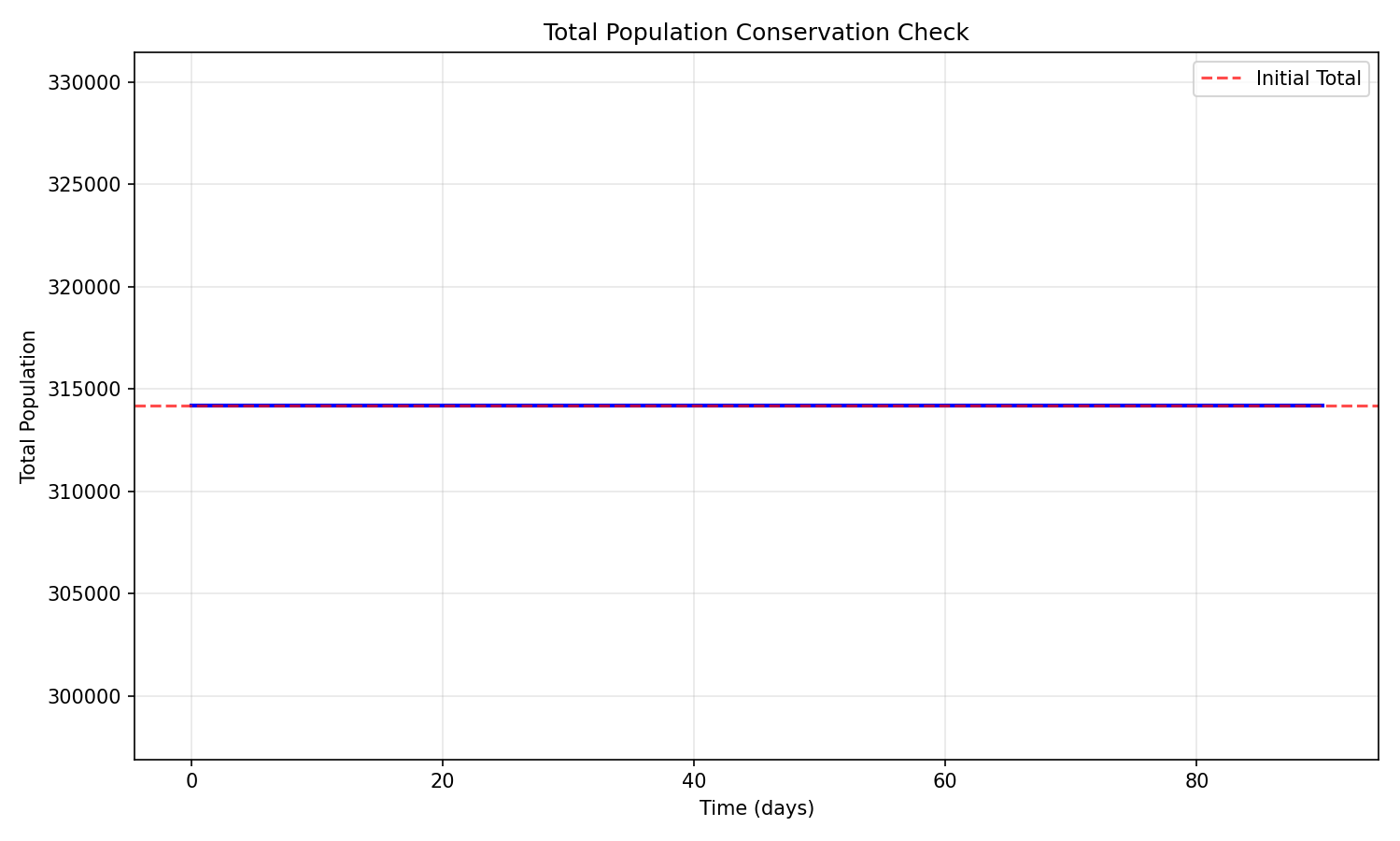}
\captionof{figure}{Total population conservation check under the proposed model.}
\label{fig:total_population}
\end{minipage}

\subsubsection{\texorpdfstring{Equilibrium and the Effect of Spatial Dimensionality in Model~\eqref{eq:1.4}--\eqref{eq:1.5}}{Equilibrium and the Effect of Spatial Dimensionality in Model 1.4--1.5}}

Let us first consider the first equation of system~(1.4). The dimension of $\partial S$ coincides with the dimension of $S$, while the dimension of $\partial t$ coincides with the dimension of time $T$. Hence, from a dimensional standpoint,
\[
\dim\!\left(\frac{\partial S}{\partial t}\right)=\frac{[S]}{[T]}.
\]

From equilibrium considerations, the dimensions of the left- and right-hand sides of the equation must coincide; therefore, we examine the right-hand side. Since the three compartments $S$, $I$, and $R$ are defined within the total population size $N$, their dimensions must be identical. Consequently,
\[
\dim\!\left(\frac{SI}{N}\right)=\frac{[S][I]}{[N]}=[S].
\]

In this case, the dimension of the parameter $\beta$ must be:
\[
\dim \beta = \frac{1}{[T]}.
\]

For the Laplacian term in the right-hand side of the first equation, it is known that it represents the sum of second partial derivatives. Hence, we have:
\[
\dim \Delta(S) = \frac{[S]}{L^2}.
\]

Therefore, the dimension of the last term must be:
\[
\dim D_S = \frac{L^2}{[T]}.
\]

For the equations for $I$ and $R$, analogous reasoning and computational steps apply.

We now determine the equilibrium (steady) states of the system, i.e., the states for which the right-hand side vanishes. Hence,
\begin{equation}
 \frac{\partial S}{\partial t}
=\frac{\partial I}{\partial t}
=\frac{\partial R}{\partial t}
=0.
\tag{3.22}\label{eq:3.22}
\end{equation}

That is,
\[
-\beta\,\frac{S_e I_e}{N_e}+D_S\Delta(S_e)=0,
\]
\[
\beta\,\frac{S_e I_e}{N_e}-\gamma I_e+D_I\Delta(I_e)=0,
\]
\[
\gamma I_e+D_R\Delta(R_e)=0.
\]
Due to condition~\eqref{eq:3.22} one has 
...
We already know that, in the equilibrium regime,
\[
\Delta(S_e)=0,\qquad \Delta(I_e)=0,\qquad \Delta(R_e)=0.
\]
Then, from the first equation we obtain $S_e I_e=0$, and therefore either $S_e=0$ or $I_e=0$.
However, for the second equation we must set $I_e=0$; then the third equation also yields the
corresponding equilibrium. In this case (with $S_e\neq 0$), the values $S_e$ and $R_e$ may be arbitrary
constants, and the system can be written in the form
\[
i=(I-I_e),\qquad s=(S_e-S),\qquad r=(R_e-R).
\]

Then our equation for perturbations $s,i$ and $r$ takes a form:
\begin{equation}
\begin{aligned}
\frac{\partial s}{\partial t} &= -\beta\,\frac{s}{N} -\beta\,\frac{I}{N}+ D_S\,\Delta s,\\
\frac{\partial i}{\partial t} &= \beta\,\frac{s}{N} +\beta\,\frac{I}{N} - \gamma\,I + D_I\,\Delta I,\\
\frac{\partial r}{\partial t} &= \gamma\,I + D_R\,\Delta r,\\
S+I+R &= N>0,\\
\Delta(\cdot) &= \frac{\partial^2(\cdot)}{\partial x^2}+\frac{\partial^2(\cdot)}{\partial y^2}.
\end{aligned}
\tag{1.4.1}
\end{equation}

\subsubsection{\texorpdfstring{Numerical Discretization and Solution Procedure for Modle~\eqref{eq:1.4}--\eqref{eq:1.5}}{Numerical Discretization and Solution Procedure for Eqs. (1.4)--(1.5)}}\label{subsubsec:sir-no-chemotaxis-numerics}

First, we assume that the initial outbreak location $(x_0,y_0)$ is situated at the center of the square
domain, and that the first infected individual is also placed at this point. We then solve
Eqs.~(4)--(5) using an implicit finite-difference scheme. For the Laplace operator, we apply the standard
second-order spatial discretization
\begin{equation}
\Delta(\cdot)\approx
\begin{aligned}
&\frac{(\cdot)^{n}_{i+1,j}-2(\cdot)^{n}_{i,j}+(\cdot)^{n}_{i-1,j}}{\Delta x^{2}}
+\frac{(\cdot)^{n}_{i,j+1}-2(\cdot)^{n}_{i,j}+(\cdot)^{n}_{i,j-1}}{\Delta y^{2}}.
\end{aligned}
\tag{2.19}
\end{equation}

For the time discretization, we simultaneously employ a fully implicit scheme and obtain
\begin{equation}
\frac{(\cdot)^{n+1}_{i,j}-(\cdot)^{n}_{i,j}}{\Delta t}
=
f\!\left((\cdot)^{n+1}_{i,j}\right)
+
D\,\Delta\!\left((\cdot)^{n+1}_{i,j}\right).
\tag{2.20}
\end{equation}

For different choices of the function $f$, which represent the reaction terms corresponding to the three
equations of the SIR model (excluding the diffusion terms), nonlinearities arise after discretization.

We employ a Picard (fixed-point) iteration, which leads to the following linear systems at step $k=0,1,2,\ldots$:
\begin{equation}
\begin{aligned}
\left[
\frac{I}{\Delta t}
+\beta\,\mathrm{diag}\!\left(\frac{I_k}{N}\right)
-D_S L_S
\right] S_{k+1}
=
\frac{S^{n}}{\Delta t},\\[4pt]
\left[
\frac{I}{\Delta t}
+\gamma I
-\beta\,\mathrm{diag}\!\left(\frac{S_{k+1}}{N}\right)
-D_I L_I
\right] I_{k+1}
=
\frac{I^{n}}{\Delta t},\\[4pt]
\left[
\frac{I}{\Delta t}
-D_R L_R
\right] R_{k+1}
=
\frac{R^{n}}{\Delta t}
+\gamma I_{k+1}.
\end{aligned}
\tag{2.21}
\end{equation}

To better illustrate the results of our experiment, we adopt the following approach: it is more informative to plot level sets (contour lines). (Here we display only the contours for (S).) In our case, these level sets are simple circles. As can be seen, the contour lines become progressively more regular at each time step. 

\textbf{At the initial time (t=0):}
\clearpage

\noindent\begin{minipage}{\linewidth}
\centering
\includegraphics[width=\linewidth]{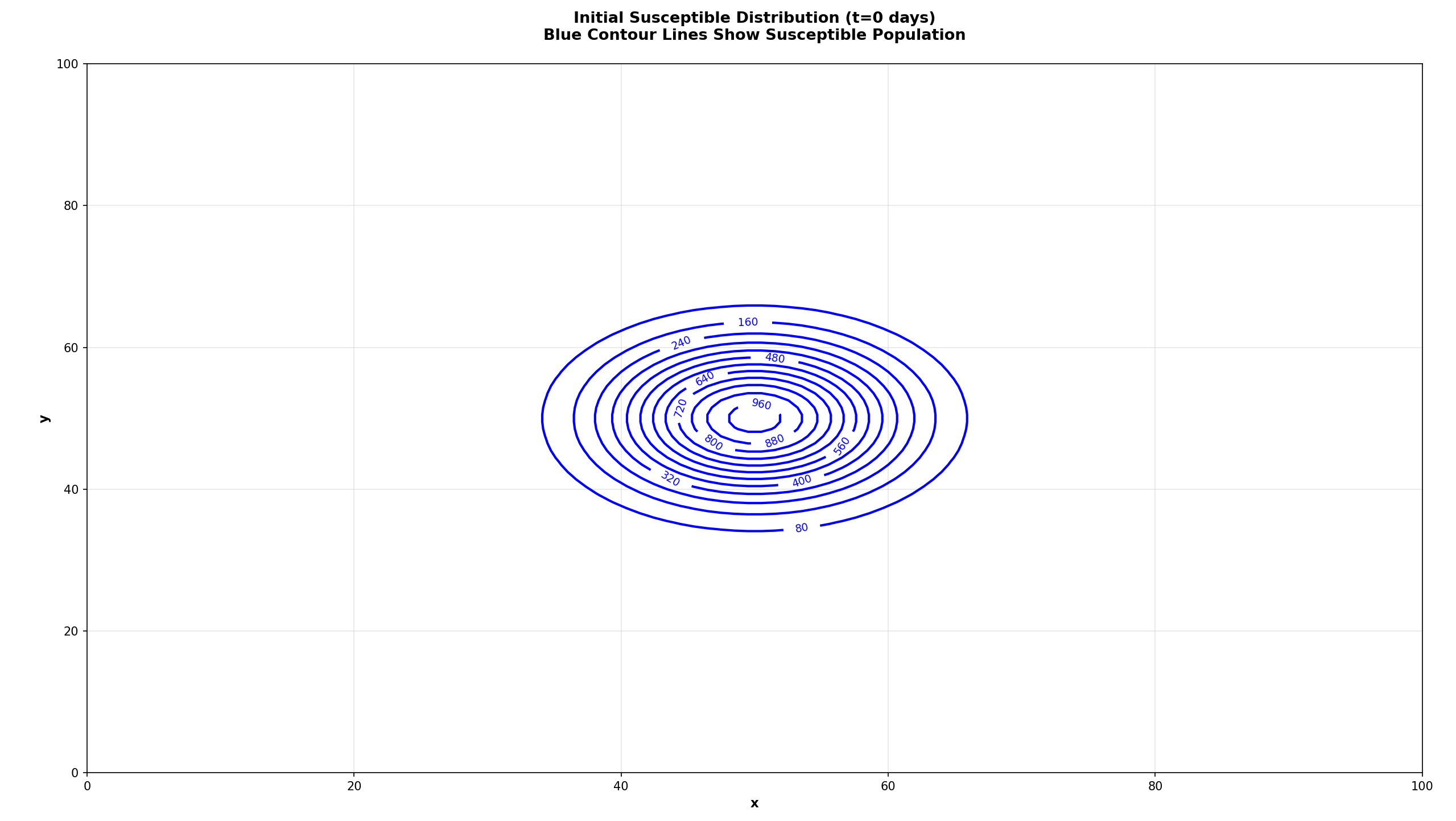}
\captionof{figure}{Initial susceptible distribution at $t=0$ (contour lines of $S$) for the reaction--diffusion SIR model without chemotaxis.}
\label{fig:pde_sir_no_chem_t0}
\end{minipage}

\textbf{At the 29th day (t=28.7):}

\noindent\begin{minipage}{\linewidth}
\centering
\includegraphics[width=\columnwidth]{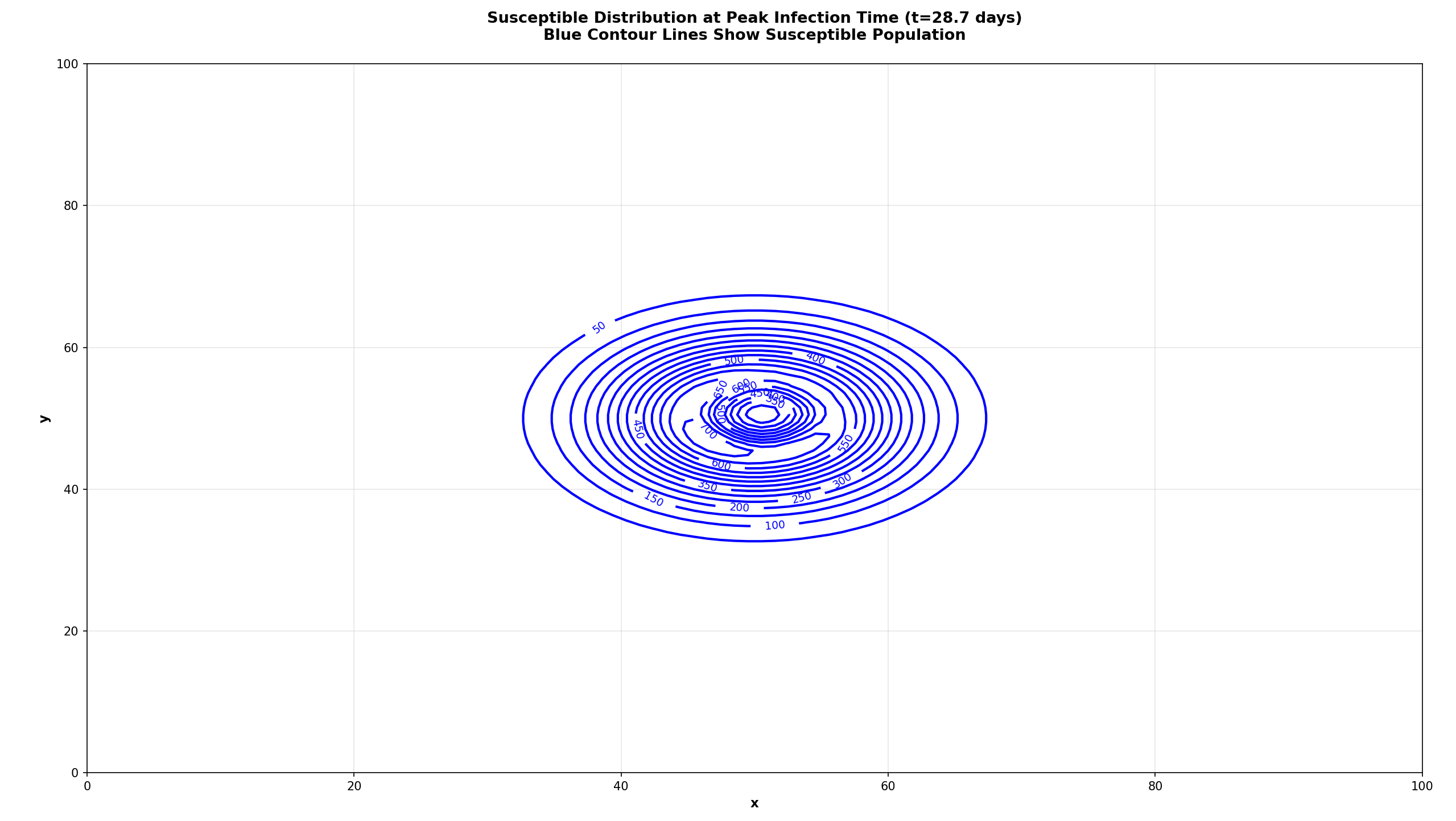}
\captionof{figure}{Susceptible distribution at peak infection time (contour lines of $S$) for the reaction--diffusion SIR model without chemotaxis.}
\label{fig:pde_sir_no_chem_t28.7}
\end{minipage}

\textbf{At the end of the experimental period (t=90):}

\noindent\begin{minipage}{\linewidth}
\centering
\includegraphics[width=\columnwidth]{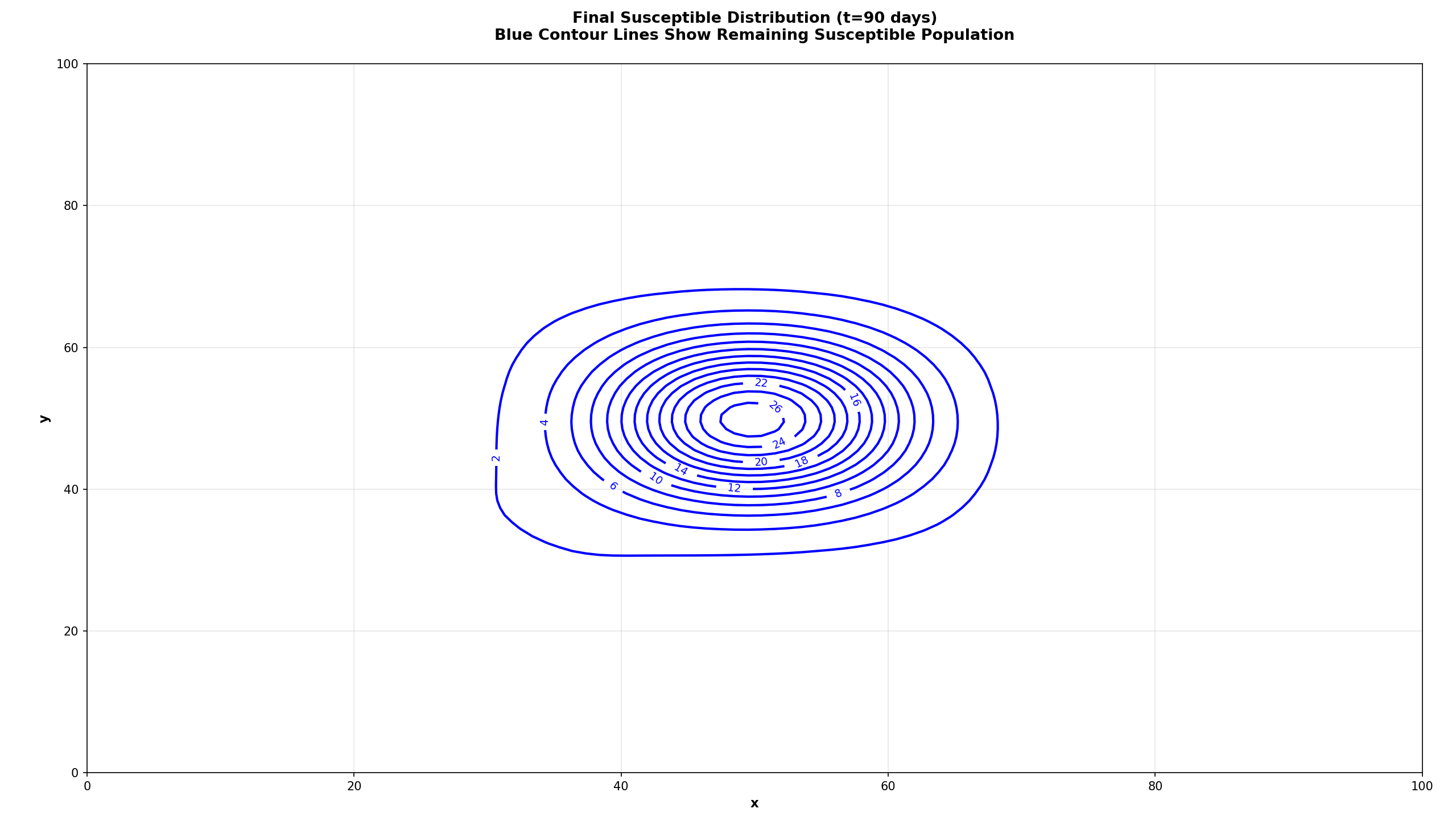}
\captionof{figure}{Final susceptible distribution at $t=90$ (contour lines of $S$) for the reaction--diffusion SIR model without chemotaxis.}
\label{fig:pde_sir_no_chem_t90}
\end{minipage}

\subsection{SIR System with a Chemotaxis Term}
\subsubsection{\texorpdfstring{Model~\eqref{eq:1.6}--\eqref{eq:1.7}}{Model 1.6--1.7}}

These figures indicate that we have already obtained satisfactory results. In the subsequent analysis of the SIR model, we consider the following aspect: the influence of the spatial location of the infection source. In particular, we study the case in which the initial distribution of infected individuals is not concentrated at the center of the domain, but rather at an arbitrary point in space. At the same time, and more importantly, once we learn that people around us are infected, we instinctively attempt to keep a certain distance from them; this behavior should be reflected by an additional term in the first equation. In this way, we systematize these ideas and obtain models~\eqref{eq:1.6}--\eqref{eq:1.7}.

\begin{equation}
\begin{alignedat}{1}
&\tau_S\,\frac{\partial S}{\partial t}
= -\beta\,\frac{S I}{N}
+ D_S\,\Delta S
+\chi\,\nabla\!\left(\frac{S\,\nabla I}{I}\right),\\
&\tau_I\,\frac{\partial I}{\partial t}
= \beta\,\frac{S I}{N}
-\gamma I
+ D_I\,\Delta I,\\
&\tau_R\,\frac{\partial R}{\partial t}
= \gamma I
+ D_R\,\Delta R,\\
&S+I+R=N>0,\\[2pt]
&\Delta(\cdot)=\frac{\partial^2(\cdot)}{\partial x^2}
+\frac{\partial^2(\cdot)}{\partial y^2},\\
&\nabla(\cdot)=\left(\frac{\partial(\cdot)}{\partial x},
\frac{\partial(\cdot)}{\partial y}\right),\\
&\nabla\!\left(\frac{S\,\nabla I}{I}\right)
=\frac{\partial}{\partial x}\!\left(\frac{S}{I}\frac{\partial I}{\partial x}\right)
+\frac{\partial}{\partial y}\!\left(\frac{S}{I}\frac{\partial I}{\partial y}\right),\\
&\tau_{(\cdot)}=1 \quad \text{is assumed for the moment}.
\end{alignedat}
\tag{1.6}\label{eq:1.6}
\end{equation}
Where $\chi$ denotes the fear level.
\begin{equation}
\begin{alignedat}{1}
&\left.\frac{\partial S}{\partial n}\right|_{\Gamma}
=\left.\frac{\partial I}{\partial n}\right|_{\Gamma}
=\left.\frac{\partial R}{\partial n}\right|_{\Gamma}
=0,\\[2pt]
&S(x,y,0)=S_0(x,y)
=A\exp\!\Bigl[-B\bigl((x-x_0)^2+(y-y_0)^2\bigr)\Bigr],\\[2pt]
&I(x,y,0)=I_0(x,y)
=E\exp\!\Bigl[-F\bigl((x-x_1)^2+(y-y_1)^2\bigr)\Bigr],\\
&R(x,y,0)=0,\\[2pt]
&\frac{\partial(\cdot)}{\partial n}
=\cos(\vec n,x)\frac{\partial(\cdot)}{\partial x}
+\cos(\vec n,y)\frac{\partial(\cdot)}{\partial y}.
\end{alignedat}
\tag{1.7}\label{eq:1.7}
\end{equation}

\subsubsection{\texorpdfstring{Numerical Discretization and Solution Procedure of Model~\eqref{eq:1.6}--\eqref{eq:1.7}}{Numerical Discretization and Solution Procedure of Model 1.6--1.7}}

We discretize the spatial domain $\Omega=[0,L_x]\times[0,L_y]$ by a uniform grid
$x_i=i\Delta x$ ($i=0,\dots,N_x-1$), $y_j=j\Delta y$ ($j=0,\dots,N_y-1$), and set $t^n=n\Delta t$.
Homogeneous Neumann (zero--flux) boundary conditions are enforced by mirror treatment at the boundary,
consistent with a discrete Neumann Laplacian.

\paragraph{Discrete Laplacian.}
Let $u_{i,j}^n$ denote a generic grid function. On interior nodes we use the standard second--order
five--point stencil
\begin{equation}
(\Delta_h u)_{i,j}=
\begin{aligned}
\frac{u_{i+1,j}-2u_{i,j}+u_{i-1,j}}{\Delta x^2}
+\frac{u_{i,j+1}-2u_{i,j}+u_{i,j-1}}{\Delta y^2},
\end{aligned}
\tag{2.22}
\end{equation}
and apply the Neumann condition by reflection at the boundary. The corresponding sparse 2D Laplacian
matrix is assembled once as a Kronecker sum of 1D Neumann Laplacians.

\paragraph{Chemotaxis (fear/avoidance) flux term.}
The additional transport term is written in conservative form as
$Q(S,I)=\nabla\cdot\!\big(S\nabla I/I\big)$.
To avoid division by small $I$, we regularize $I$ by
$I_{\mathrm{safe}}=\max(I,\varepsilon)$ and use $\nabla I/I=\nabla(\ln I_{\mathrm{safe}})$.
Hence,
\begin{equation}
Q(S,I)=\nabla\cdot\!\Big(S\,\frac{\nabla I}{I}\Big)
      =\nabla\cdot\!\big(S\nabla(\ln I_{\mathrm{safe}})\big).
\tag{2.23}\label{eq:2.23}
\end{equation}

Let $w=\ln I_{\mathrm{safe}}$. We compute $\nabla w$ by central differences and then approximate the
divergence of the flux $S\nabla w$ (again by central differences):
\begin{equation}
\begin{aligned}
(w_x)_{i,j} &\approx \frac{w_{i+1,j}-w_{i-1,j}}{2\Delta x},\\
(w_y)_{i,j} &\approx \frac{w_{i,j+1}-w_{i,j-1}}{2\Delta y},\\[4pt]
(F_x)_{i,j} &= S_{i,j}(w_x)_{i,j},\\
(F_y)_{i,j} &= S_{i,j}(w_y)_{i,j},\\[4pt]
Q_{i,j} &\approx
\frac{(F_x)_{i+1,j}-(F_x)_{i-1,j}}{2\Delta x} +
\frac{(F_y)_{i,j+1}-(F_y)_{i,j-1}}{2\Delta y}.
\end{aligned}
\tag{2.24}\label{eq:2.24}
\end{equation}
In the implementation, the same mirror padding is used at the boundary to remain compatible with the
zero--flux condition.

\paragraph{Time stepping (semi-implicit IMEX Euler).}
Denote $N^n_{i,j}=S^n_{i,j}+I^n_{i,j}+R^n_{i,j}$ and define the (explicit) infection term
\begin{equation}
f^n_{i,j}=\beta\,\frac{S^n_{i,j}I^n_{i,j}}{\max(N^n_{i,j},10^{-10})}.
\tag{2.25}
\end{equation}
We employ a first--order IMEX Euler scheme: diffusion and linear recovery are treated implicitly,
whereas the nonlinear infection term and the chemotaxis term are evaluated explicitly at time level
$n$. The update reads
\begin{equation}
\Big(\frac{\tau_S}{\Delta t}-D_S\Delta_h\Big) S^{n+1}
=
\frac{\tau_S}{\Delta t}S^{n} - f^{n} + \chi\,Q^{n},
\tag{2.26}\label{eq:2.26}
\end{equation}
\begin{equation}
\Big(\frac{\tau_I}{\Delta t}+\gamma-D_I\Delta_h\Big) I^{n+1}
=
\frac{\tau_I}{\Delta t}I^{n} + f^{n},
\tag{2.27}\label{eq:2.27}
\end{equation}
\begin{equation}
\Big(\frac{\tau_R}{\Delta t}-D_R\Delta_h\Big) R^{n+1}
=
\frac{\tau_R}{\Delta t}R^{n} + \gamma\,I^{n+1}.
\tag{2.28}\label{eq:2.28}
\end{equation}
Here $Q^n=Q(S^n,I^n)$ is computed via~\eqref{eq:2.23}--\eqref{eq:2.24}.

\paragraph{Sparse linear solves and stabilization.}
After flattening the $N_x\times N_y$ fields into vectors,~\eqref{eq:2.26}--\eqref{eq:2.28} become sparse linear systems with time--independent matrices
\begin{equation}
\begin{aligned}
A_S &= \Big(\frac{\tau_S}{\Delta t}\Big)I - D_S L_{2D},\\
A_I &= \Big(\frac{\tau_I}{\Delta t}+\gamma\Big)I - D_I L_{2D},\\
A_R &= \Big(\frac{\tau_R}{\Delta t}\Big)I - D_R L_{2D}.
\end{aligned}
\tag{2.29}
\end{equation}
which are solved at each time step using a sparse direct solver. To preserve physical admissibility,
we enforce nonnegativity by projection:
\begin{equation}
\begin{aligned}
S^{n+1} &\leftarrow \max(S^{n+1},0),\\
I^{n+1} &\leftarrow \max(I^{n+1},0),\\
R^{n+1} &\leftarrow \max(R^{n+1},0).
\end{aligned}
\tag{2.30}
\end{equation}
Finally, an optional uniform rescaling is applied to reduce accumulated mass drift:
\begin{equation}
M^n=\sum_{i,j}\big(S^n_{i,j}+I^n_{i,j}+R^n_{i,j}\big)\,\Delta x\,\Delta y.
\tag{2.31}
\end{equation}
\begin{equation}
\alpha=\frac{M^0}{M^{n+1}}.
\tag{2.32}
\end{equation}
\begin{equation}
(S^{n+1},I^{n+1},R^{n+1})\leftarrow \alpha\,(S^{n+1},I^{n+1},R^{n+1}).
\tag{2.33}
\end{equation}

In this case, we assume that the initial infection source is located in the northwestern part of the city; accordingly, the evolution of the susceptible population over the time interval (0!-!90) days is as follows. 

\textbf{At the initial time (t=0):}

\noindent\begin{minipage}{\linewidth}
\centering
\includegraphics[width=\columnwidth]{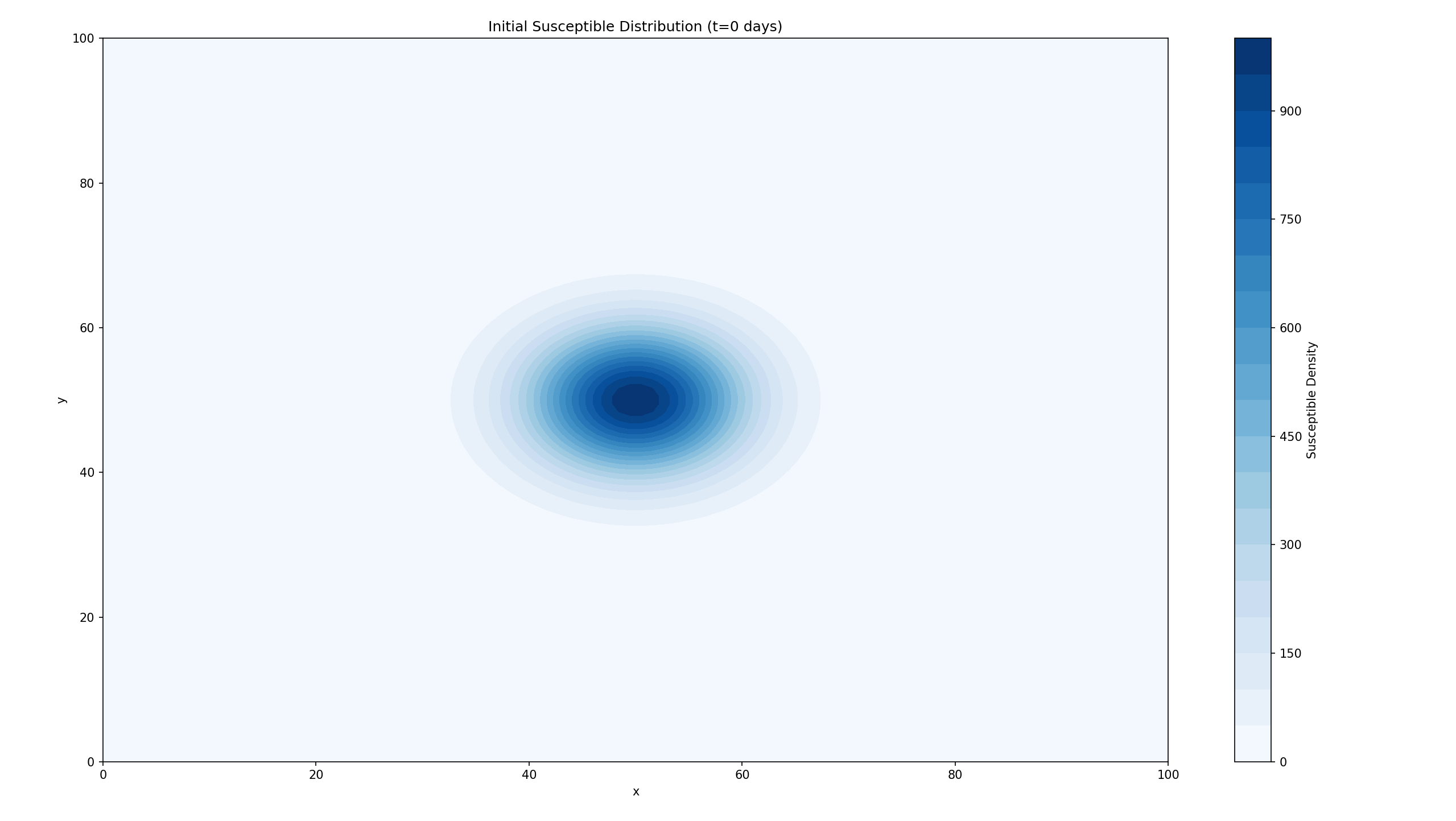}
\captionof{figure}{Initial susceptible distribution at $t=0$ (contour lines of $S$) for the reaction--diffusion SIR model with chemotaxis.}
\label{fig:pde_sir_wth_chem_t0}
\end{minipage}

\textbf{At the 40th day (t=39.5):}

\noindent\begin{minipage}{\linewidth}
\centering
\includegraphics[width=\columnwidth]{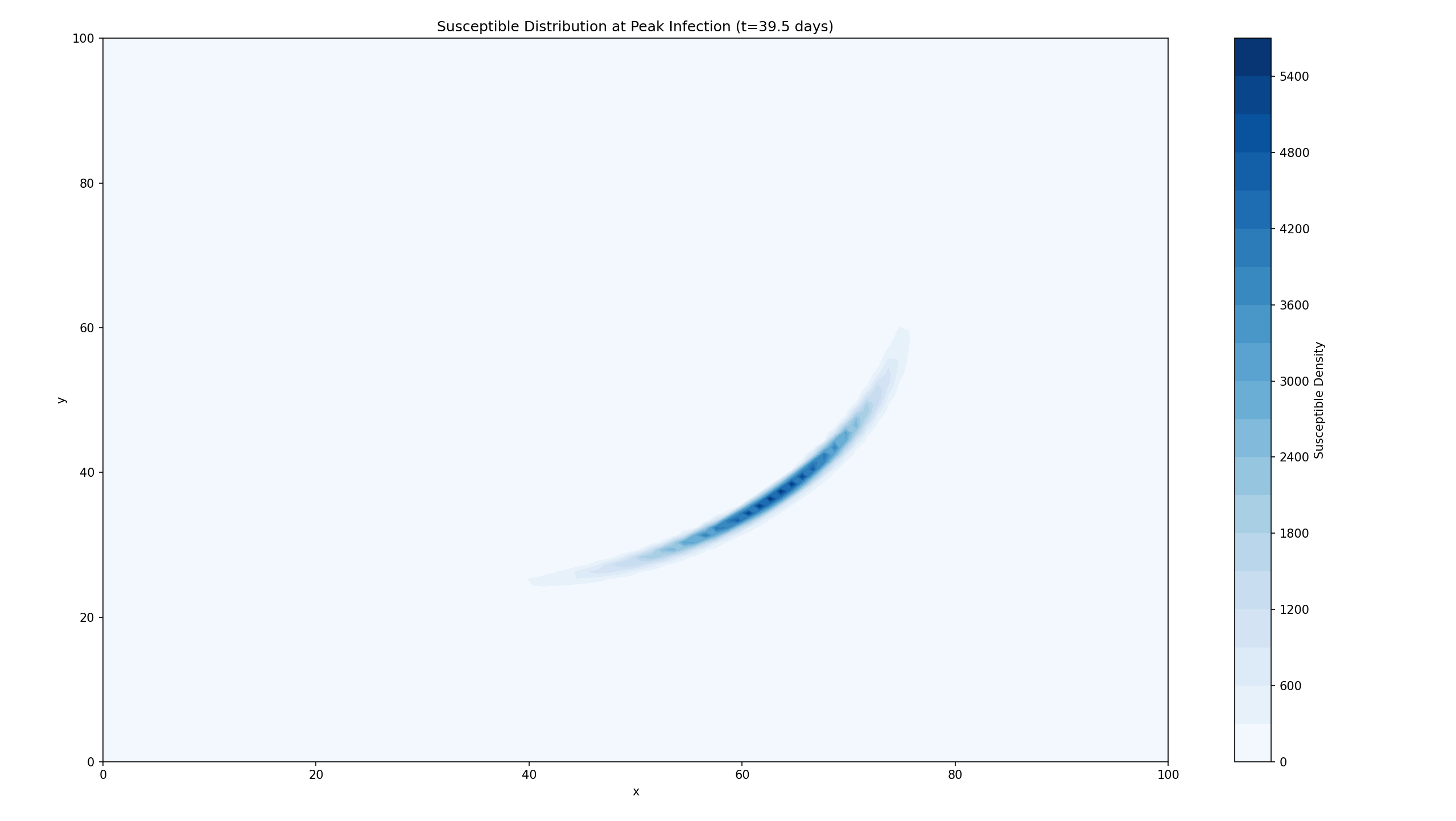}
\captionof{figure}{Susceptible distribution at peak infection time (contour lines of $S$) for the reaction--diffusion SIR model with chemotaxis.}
\label{fig:pde_sir_with_chem_t40}
\end{minipage}

\textbf{At the end of the experimental period (t=90):}

\noindent\begin{minipage}{\linewidth}
\centering
\includegraphics[width=\columnwidth]{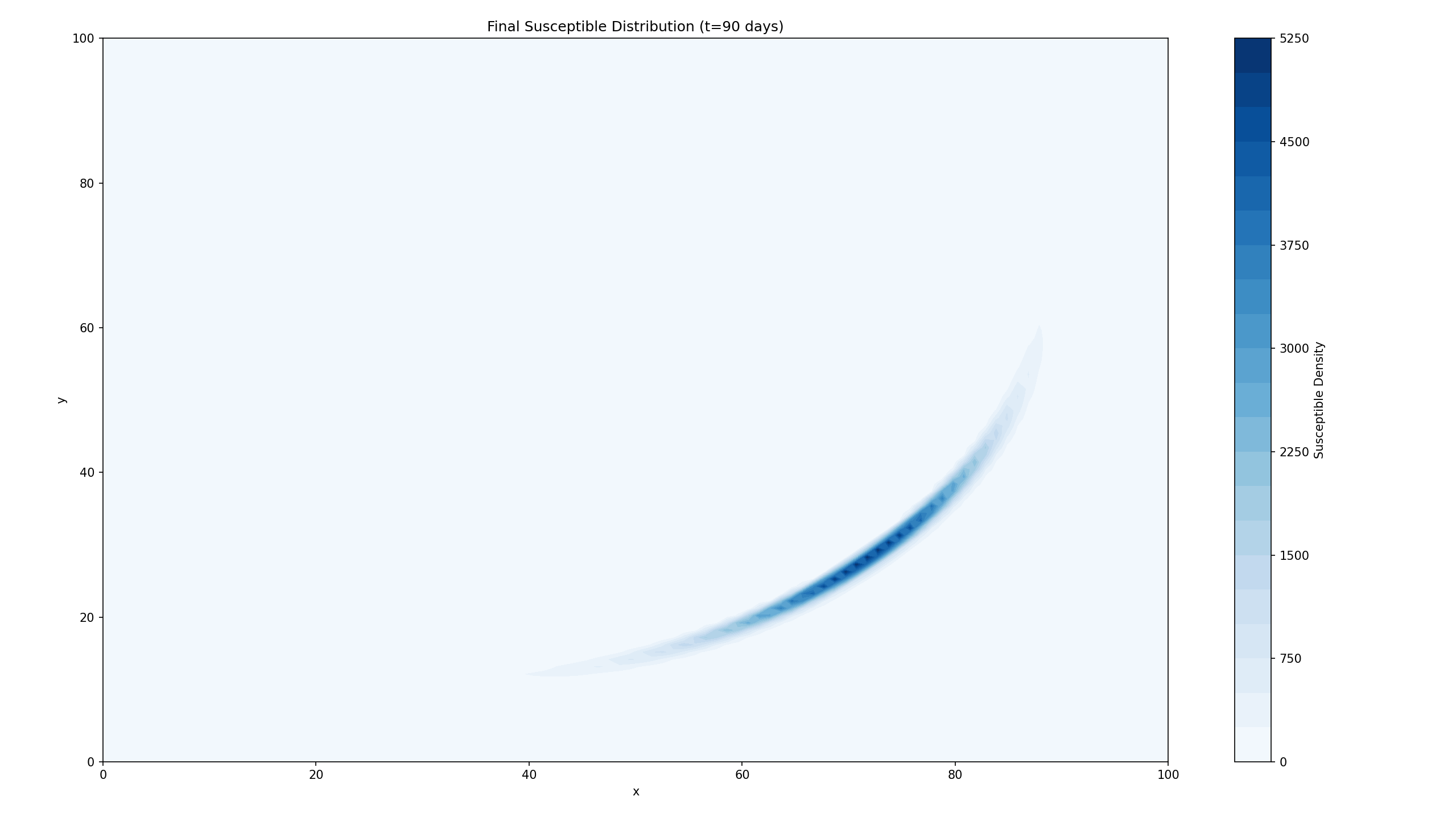}
\captionof{figure}{Final susceptible distribution at $t=90$ (contour lines of $S$) for the reaction--diffusion SIR model with chemotaxis.}
\label{fig:pde_sir_with_chem_t90}
\end{minipage}

\subsubsection{\texorpdfstring{Application of the Einstein Paradigm to Model~\eqref{eq:1.6}--\eqref{eq:1.7}}{Application of the Einstein Paradigm to Model 1.6--1.7}}

We continue modeling the transport line for this model, considering an arbitrary trajectory of a bus moving between two stops, as an example of such a curve:

\noindent\begin{minipage}{\linewidth}
\centering
\includegraphics[width=0.7\linewidth]{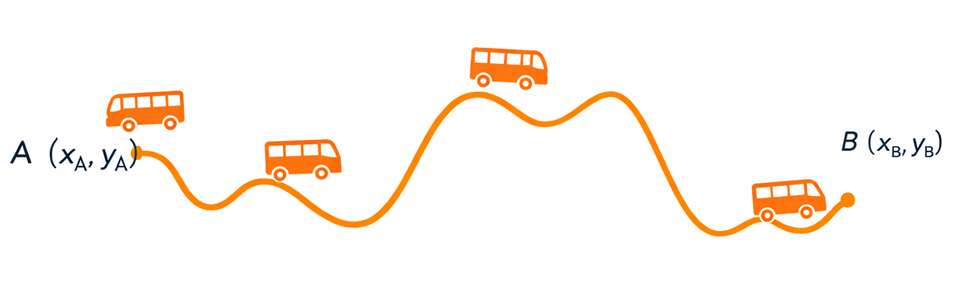}
\captionof{figure}{Radom bus route between two stops.}
\label{fig:bus_route}
\end{minipage}

By analogy with~(3.1), we generalize it to the two-dimensional case and obtain
\begin{equation}
\begin{aligned}
S(x,y,t+\tau)
&=\int_{-\infty}^{+\infty}\!\!\int_{-\infty}^{+\infty}
\Big(S(x+\zeta_1,y+\zeta_2,t)-S(x,y,t)\Big)\,
\varphi(x,y,t,\zeta_1,\zeta_2)\,d\zeta_1\,d\zeta_2 \\
&\quad + S(x,y,t)\int_{-\infty}^{+\infty}\!\!\int_{-\infty}^{+\infty}
\varphi(x,y,t,\zeta_1,\zeta_2)\,d\zeta_1\,d\zeta_2
-\beta\,\frac{S I}{N}.
\end{aligned}
\tag{4.5}
\end{equation}

Expanding the first difference term into a Taylor series (up to second order) yields
\begin{equation}
\begin{aligned}
&S(x+\zeta_1,y+\zeta_2,t)-S(x,y,t)
=\zeta_1\frac{\partial S}{\partial x}+\zeta_2\frac{\partial S}{\partial y}
+\frac{\zeta_1^2}{2}\frac{\partial^2 S}{\partial x^2}
+\frac{\zeta_2^2}{2}\frac{\partial^2 S}{\partial y^2}
+\frac{\zeta_1\zeta_2}{2}\frac{\partial^2 S}{\partial x\,\partial y}
+\frac{\zeta_2\zeta_1}{2}\frac{\partial^2 S}{\partial y\,\partial x}
+\cdots \\
&=\frac{\partial S}{\partial x}\int_{-\infty}^{+\infty}\!\!\int_{-\infty}^{+\infty}
\zeta_1\,\varphi_S(x,y,t,\zeta_1,\zeta_2)\,d\zeta_1\,d\zeta_2
+\frac{\partial S}{\partial y}\int_{-\infty}^{+\infty}\!\!\int_{-\infty}^{+\infty}
\zeta_2\,\varphi_S(x,y,t,\zeta_1,\zeta_2)\,d\zeta_1\,d\zeta_2 \\
&\quad +\frac{\partial^2 S}{\partial x^2}\int_{-\infty}^{+\infty}\!\!\int_{-\infty}^{+\infty}
\frac{\zeta_1^2}{2}\,\varphi_S(x,y,t,\zeta_1,\zeta_2)\,d\zeta_1\,d\zeta_2
+\frac{\partial^2 S}{\partial y^2}\int_{-\infty}^{+\infty}\!\!\int_{-\infty}^{+\infty}
\frac{\zeta_2^2}{2}\,\varphi_S(x,y,t,\zeta_1,\zeta_2)\,d\zeta_1\,d\zeta_2
+\cdots .
\end{aligned}
\tag{4.6}
\end{equation}

For simplicity, we assume that
\begin{equation}
\boxed{
\int_{-\infty}^{+\infty}\!\!\int_{-\infty}^{+\infty}
\zeta_1\zeta_2\,\varphi_S(x,y,t,\zeta_1,\zeta_2)\,d\zeta_1\,d\zeta_2=0
}.
\tag{4.7}
\end{equation}

Moreover, we continue to neglect the Taylor remainder term and formally set it to zero:
\begin{equation}
...=R_2(x,y)\approx 0.
\tag{4.8}
\end{equation}

Therefore, the last two contributions in~(4.6) vanish. Hence, for the equation for $S$ we can write:
\begin{equation}
\boxed{
\begin{aligned}
&S(x,y,t+\tau)-S(x,y,t)
=S(x,y,t)\,\frac{\partial}{\partial x}
\int_{-\infty}^{+\infty}\!\!\int_{-\infty}^{+\infty}
\zeta_1\,\varphi_S(x,y,t,\zeta_1,\zeta_2)\,d\zeta_1\,d\zeta_2 \\
&\quad +S(x,y,t)\,\frac{\partial}{\partial y}
\int_{-\infty}^{+\infty}\!\!\int_{-\infty}^{+\infty}
\zeta_2\,\varphi_S(x,y,t,\zeta_1,\zeta_2)\,d\zeta_1\,d\zeta_2 \\
&\quad +\tau_S\,\overline{F}_S\!\big(S, I(x,y,t)\big)
+\int_{-\infty}^{+\infty}\!\!\int_{-\infty}^{+\infty}
\Big(S(x+\zeta_1,y+\zeta_2,t)-S(x,y,t)\Big)\,
\varphi_S(x,y,t,\zeta_1,\zeta_2)\,d\zeta_1\,d\zeta_2 .
\end{aligned}
}
\tag{4.9}
\end{equation}

We expand the last integral in~(4.9); as a result, we obtain an expression containing second--order derivatives:
\begin{equation}
\boxed{
\begin{aligned}
&S(x,y,t+\tau)-S(x,y,t)
=\textcolor{blue}{S(x,y,t)\,\frac{\partial}{\partial x}\int_{-\infty}^{+\infty}\!\!\int_{-\infty}^{+\infty}
\zeta_1\,\varphi_S(x,y,t,\zeta_1,\zeta_2)\,d\zeta_1\,d\zeta_2}
\;+\\
&\quad+\textcolor{blue}{S(x,y,t)\,\frac{\partial}{\partial y}\int_{-\infty}^{+\infty}\!\!\int_{-\infty}^{+\infty}
\zeta_2\,\varphi_S(x,y,t,\zeta_1,\zeta_2)\,d\zeta_1\,d\zeta_2}+\\
&\quad+\textcolor{blue}{\int_{-\infty}^{+\infty}\!\!\int_{-\infty}^{+\infty}
\zeta_2\,\varphi_S(x,y,t,\zeta_1,\zeta_2)\,d\zeta_1\,d\zeta_2\;\frac{\partial S}{\partial y}}
\textcolor{blue}{+\int_{-\infty}^{+\infty}\!\!\int_{-\infty}^{+\infty}
\zeta_1\,\varphi_S(x,y,t,\zeta_1,\zeta_2)\,d\zeta_1\,d\zeta_2\;\frac{\partial S}{\partial x}}\\
&\quad+\frac{\partial^2 S}{\partial x^2}\int_{-\infty}^{+\infty}\!\!\int_{-\infty}^{+\infty}
\frac{\zeta_1^{2}}{2}\,\varphi_S(x,y,t,\zeta_1,\zeta_2)\,d\zeta_1\,d\zeta_2
+\frac{\partial^2 S}{\partial y^2}\int_{-\infty}^{+\infty}\!\!\int_{-\infty}^{+\infty}
\frac{\zeta_2^{2}}{2}\,\varphi_S(x,y,t,\zeta_1,\zeta_2)\,d\zeta_1\,d\zeta_2\\
&\quad+\tau_S\,\overline{F}_S\!\big(S, I(x,y,t)\big).
\end{aligned}
}
\tag{4.10}
\end{equation}

We extend Hypothesis~2.1 to the two--dimensional case by formulating Hypothesis~3.1, and obtain:
\begin{hypothesis}
\begin{equation}
\boxed{
\begin{aligned}
\frac{\partial}{\partial x}\int_{-\infty}^{+\infty}\!\!\int_{-\infty}^{+\infty}
\zeta_1\,\varphi_S(x,y,t,\zeta_1,\zeta_2)\,d\zeta_1\,d\zeta_2
&=-\chi\,\frac{1}{I}\,\frac{\partial I}{\partial x},\\[4pt]
\frac{\partial}{\partial y}\int_{-\infty}^{+\infty}\!\!\int_{-\infty}^{+\infty}
\zeta_2\,\varphi_S(x,y,t,\zeta_1,\zeta_2)\,d\zeta_1\,d\zeta_2
&=-\chi\,\frac{1}{I}\,\frac{\partial I}{\partial y}.
\end{aligned}
}
\tag{4.11}
\end{equation}
\end{hypothesis}

Thus, in~(4.10) we group the terms that do not contain second--order derivatives (highlighted in blue) by first--order partial derivatives, and obtain the chemotactic term of our model:
\begin{equation}
\boxed{
\begin{aligned}
&\textcolor{blue}{
S(x,y,t)\,\frac{\partial}{\partial x}\!\int_{-\infty}^{+\infty}\!\!\int_{-\infty}^{+\infty}
\zeta_1\,\varphi_S(x,y,t,\zeta_1,\zeta_2)\,d\zeta_1\,d\zeta_2} \\
&\textcolor{blue}{\quad+S(x,y,t)\,\frac{\partial}{\partial y}\!\int_{-\infty}^{+\infty}\!\!\int_{-\infty}^{+\infty}
\zeta_2\,\varphi_S(x,y,t,\zeta_1,\zeta_2)\,d\zeta_1\,d\zeta_2 }\\
&\quad\textcolor{blue}{\;+\int_{-\infty}^{+\infty}\!\!\int_{-\infty}^{+\infty}
\zeta_2\,\varphi_S(x,y,t,\zeta_1,\zeta_2)\,d\zeta_1\,d\zeta_2\;\frac{\partial S}{\partial y}
+\int_{-\infty}^{+\infty}\!\!\int_{-\infty}^{+\infty}
\zeta_1\,\varphi_S(x,y,t,\zeta_1,\zeta_2)\,d\zeta_1\,d\zeta_2\;\frac{\partial S}{\partial x}}\\
&\qquad=\chi\,\frac{\partial}{\partial x}\!\left(S\,\frac{1}{I}\frac{\partial I}{\partial x}\right)
+\chi\,\frac{\partial}{\partial y}\!\left(S\,\frac{1}{I}\frac{\partial I}{\partial y}\right).
\end{aligned}
}
\tag{4.12}
\end{equation}

Since the chemotactic term is absent in the equations for $I$ and $R$, we directly extend their one--dimensional formulation within the Einstein paradigm to the two--dimensional case, i.e.,
\begin{equation}
\boxed{
\begin{aligned}
&I(x,y,t+\tau)-I(x,y,t)
=\tau_I\,\overline{F}_I\!\big(S,I(x,y,t)\big) \\
&\quad+\int_{-\infty}^{+\infty}\!\!\int_{-\infty}^{+\infty}
\Big(I(x+\zeta_1,y+\zeta_2,t)-I(x,y,t)\Big)\,
\varphi_I(x,y,t,\zeta_1,\zeta_2)\,d\zeta_1\,d\zeta_2 .
\end{aligned}
}
\tag{4.13}
\end{equation}

\begin{equation}
\boxed{
\begin{aligned}
&R(x,y,t+\tau)-R(x,y,t)
=\tau_R\,\overline{F}_R\!\big(I,R(x,y,t)\big) \\
&\quad+\int_{-\infty}^{+\infty}\!\!\int_{-\infty}^{+\infty}
\Big(R(x+\zeta_1,y+\zeta_2,t)-R(x,y,t)\Big)\,
\varphi_R(x,y,t,\zeta_1,\zeta_2)\,d\zeta_1\,d\zeta_2 .
\end{aligned}
}
\tag{4.14}
\end{equation}

Consequently, based on~(4.10), (4.13)--(4.14), we successfully construct the SIR model with a chemotactic term, namely models~\eqref{eq:1.6}--\eqref{eq:1.7}.

\section*{Acknowledgments}
The author gratefully acknowledges \textbf{Akif Ibragimov} for valuable discussions and helpful guidance on the mathematical and technical aspects of this work. The author also sincerely thanks \textbf{Evgenia Echkina} for her careful supervision, helpful comments, and support in refining the central research idea. Appreciation is further extended to the Faculty of Computational Mathematics and Cybernetics, Lomonosov Moscow State University, for providing a stimulating academic environment and favorable conditions for study and research. Finally, the author is deeply grateful to his parents for their constant support, and to \textbf{Irina Zakharova} for her recommendation and encouragement, which played an important role in his academic path.

\bibliographystyle{plainnat-revised}
\bibliography{1}


\end{document}